\documentclass[a4paper,11pt]{article}
\usepackage[hmargin=2.5cm,vmargin=3cm]{geometry}
\makeatletter
\let\@fnsymbol\@arabic
\makeatother

\usepackage{amsmath,amssymb,amsfonts,amscd,amsthm}
\usepackage{ascmac} 
\usepackage[mathscr]{eucal}
\usepackage{enumerate}
\usepackage{indentfirst} 
\usepackage{graphicx,color}
\newtheorem{thm}{Theorem}[section]
\newtheorem{prop}[thm]{Proposition}
\newtheorem{cor}[thm]{Corollary}
\newtheorem{dfn}[thm]{Definition}
\newtheorem{lem}[thm]{Lemma}
\newtheorem{ex}{Example}[section]
\newtheorem{rmk}{Remark}[section]
\numberwithin{equation}{section}
\allowdisplaybreaks[4]

\title{GEOMETRY OF BUNDLE-VALUED MULTISYMPLECTIC STRUCTURES WITH LIE ALGEBROIDS}
\author{Yuji Hirota\thanks{hirota@azabu-u.ac.jp; Division of Integrated Sciences, Azabu University, Sagamihara, Kanagawa 252-5201, Japan.}
 \and 
Noriaki Ikeda\thanks{nikeda@se.ritsumei.ac.jp; Department of Mathematical Sciences, 
Ritsumeikan University,Kusatsu, Shiga 525-8577, Japan.}
}
\date{}
\begin{document}
\maketitle 

\begin{abstract}
We study multisymplectic structures taking values in vector bundles with connections from the viewpoint of the Hamiltonian symmetry.  
We introduce the notion of bundle-valued $n$-plectic structures and exhibit some properties of them. 
In addition, we define bundle-valued homotopy momentum sections for bundle-valued $n$-plectic manifolds with Lie algebroids 
to discuss momentum map theories in both cases of quaternionic K\"{a}hler manifolds and hyper-K\"{a}hler manifolds. 
Furthermore, we generalize the Marsden-Weinstein-Meyer reduction theorem for symplectic manifolds and construct two kinds of reductions of vector-valued 1-plectic manifolds. 
\end{abstract}
\tableofcontents 
\section{Introduction}

Multisymplectic geometry is a higher analog of symplectic geometry, which originates from the study of classical field theory 
(see \cite{Kfin73, GIMM98}, for instance, \cite{RW19} for recent developments).  
It has been developing with the aim to establish suitable geometric models describing the Hamiltonian formalisms. 
In a modern context, multisymplectic manifolds, which have principal roles in multisymplectic geometry appears in sigma models \cite{HIhom22}. 
Multisymplectic manifolds can be found in the context not only of physics but also of mathematics. 
For example, an oriented manifold is a multisymplectic manifold by the volume form. A hyper-K\"{a}hler manifold and a quaternionic K\"{a}hler 
manifold are thought of as both multisymplectic manifolds. Multisymplectic geometry might connect various subjects in geometry with each other 
and enable us to discuss them in one framework. Multisymplectic geometry has a great interest to study from a geometric point of view, also. 
Similarly to the case of symplectic geometry, the theories of momentum maps and the reduction are crucial subjects to study in multisymplectic geometry.  
There are various kinds of momentum maps for multisymplectic manifolds such as multimomentum maps \cite{MSmulti12,MSclos13}, 
homotopy moment maps \cite{CFRZhom16}, weak momentum maps \cite{Hnoe18} and so on. 
A reduction of a multisymplectic manifold has been recently discussed in the work by C. Blacker \cite{Bred21}. 

Recently, the authors have proposed homotopy momentum sections for (pre-)multisymplectic manifolds with a Lie algebroid in \cite{HIhom22}. 
It is generalization of both homotopy moment maps and momentum sections. 
Here, a momentum section is a section of a Lie algebroid with a vector bundle connection over a presymplectic manifold satisfying some conditions, 
which is introduced by C. Blohmann and A. Weinstein in \cite{BWham18}. 
Both momentum sections and homotopy momentum sections are inspired by physical analysis in which Lie groupoid (or Lie algebroid) symmetries naturally appear 
\cite{Alekseev:2004np, Blohmann:2010jd, Cattaneo:2000iw, Ikeda:2012pv, KotovStrobl19}. 
A homotopy momentum section might be a device to understand various momentum maps in a single framework. 

However, there is a drawback that a homotopy momentum section is unapplicable to the momentum map of 
quaternionic K\"{a}hler manifolds, thought of as $3$-plectic manifolds, introduced in the work of K. Galicki and L. B. Lawson, Jr. \cite{GLqua88}. 
In the paper, to unify further momentum maps theories, 
we propose a bundle-valued $n$-plectic structure and devise a homotopy momentum section for it, called a bundle-valued homotopy momentum section (BHMS for short). 
BHMSs consequently provide with us a unified framework to understand various momentum map theories including the momentum map for a quaternionic K\"{a}hler manifold. 
Additionally, we introduce the notion that BHMSs are compatible with Lie algebroids. We show that the BHMS for $n=1$ which satisfies the compatibility condition 
defines a Lie algebra structure on the linear subspace associated to a bundle-valued $1$-plectic structure. 
Furthermore, we exhibit two kinds of reduction theorems for vector-valued $1$-plectic manifolds with Lie algebroid symmetries by using BHMSs. 

The paper is organized as follows:~in Section 2, we define bundle-valued (pre-)$n$-plectic manifolds and exhibit the Cartan-type formulas for them. 
We also introduce the notion of pseudo-Hamiltonian differential forms and define the operator for them. 
In Section 3, we define a BHMS for a bundle-valued (pre-)$n$-plectic manifold and discuss some examples. 
Section 4 addresses the symmetry in quaternionic K\"{a}hler geometry and discuss the relation between the momentum maps of quaternionic K\"{a}hler manifolds and BHMSs. 
We describe the condition for the quaternionic momentum map to be a BHMS in terms of the operator defined in Section 2 (Theorem \ref{sec4:thm_G-L}). 
In Section 5, we define the compatibility with a Lie algebroid for a BHMS. In the case of $n=1$, we write the explicit equation for a BHMS to be compatible (Proposition \ref{sec5:prop_bracket}). 
Additionally, the linear subspace of the pseudo-Hamiltonian forms is shown to be a Lia algebra by the operator in Section 2 (Proposition \ref{sec5:prop_bracket2}). 
Section 6 address the reduction problem for vector-valued $1$-plectic manifolds. We shall construct the reduction in two cases: the first one is formulated
in terms of a BHMS transversal to the anchor map of a Lie algebroid (Theorem \ref{sec6:thm_main1}). 
The second one is done in terms of a BHMS satisfying the compatibility condition in Section 5 (Theorem \ref{sec6:thm_main2}). 
Both reductions are new ones, which can be applicable for vector valued $1$-plectic manifolds with Lie algebroid symmetries. 
\medskip 

In the paper, all manifolds and maps between them are assumed to be smooth. If $M$ is a smooth manifold, $C^{\infty}(M)$ denotes the space of all smooth functions on $M$, and 
$\mathfrak{X}(M)$ does the space of all smooth vector fields on $M$. For $k\geqq  0$, the space of all differential $k$-forms on $M$ is denoted by $\Omega^k(M)$. 
Given a smooth vector bundle $E\to M$, we denote by $\varGamma(E)$ the space of all smooth sections of $E$.

\section{Bundle-Valued $n$-plectic Structures}

Let $M$ be a smooth manifold and $E$ a vector bundle over $M$. We denote by $\Omega^k(M,E)$ the space of all $E$-valued $k$-forms on $M$ with $k\geqq 0$. When $k=0$, 
$\Omega^0(M,E)$ is just $\varGamma(E)$. 
Suppose that $E$ is equipped with a vector bundle connection $\nabla^E$. 
That is, $\nabla^E$ is a $\mathbb{R}$-linear mapping $\nabla^E:\varGamma(E)\to \Omega^1(M,E)$ which satisfies the Leibniz rule 
$\nabla^E (fs)={\rm d}f\otimes s + f\nabla^E s$ for any $f\in C^\infty(M)$ and $s\in \varGamma(E)$. 
We write $\nabla^E_Xs$ for $(\nabla^E s)(X)$, where $X\in\mathfrak{X}(M)$. 
\bigskip 

We now have a sequence 
\[
\cdots\overset{{\rm d}^E_\nabla}{\longrightarrow} \Omega^{k-1}(M,E) \overset{{\rm d}^E_\nabla}{\longrightarrow} \Omega^k(M,E)\overset{{\rm d}^E_{\nabla}}{\longrightarrow} 
\Omega^{k+1}(M,E)\overset{{\rm d}^E_\nabla}{\longrightarrow} \cdots\, , 
\]
where ${\rm d}^E_\nabla$ is a covariant exterior derivative, which is defined as
\begin{align}\label{sec2:eqn_cov ext der}
\bigl({\rm d}^E_{\nabla}\varphi\bigr)(X_1,\cdots,X_{p+1})
:= \sum_{i=1}^{p+1}&(-1)^{i-1}\nabla^E_{X_i}(\varphi(X_1,\cdots,\check{X}_i,\cdots,X_{p+1})) \notag \\
&+ \sum_{i<j}(-1)^{i+j}\varphi([X_i,\,X_j],X_1,\cdots,\check{X}_i,\cdots,\check{X}_j,\cdots,X_{p+1}) 
\end{align}
for $\varphi\in \Omega^k(M,E)$. The check $\check{X}_i$ means that the corresponding entry $X_i$ is omitted. 
We say that an $E$-valued form $\varphi$ is ${\rm d}^E_\nabla$-closed if ${\rm d}^E_\nabla\varphi = 0$. 

\begin{dfn}
Let $n\geqq 1$. A ${\rm d}^E_{\nabla}$-closed $E$-valued $(n+1)$-form $\omega\in \Omega^{n+1}(M,E)$ on $M$ is called an $E$-valued $n$-plectic form {\rm (}or structure{\rm )} 
if it is non-degenerate in the sense that, at each point $x\in M$, 
the induced map 
\begin{equation}
 \omega^{\flat}_x : T_xM \longrightarrow \bigwedge^{n}T^*_xM\otimes_{\mathbb{R}} E_x,\quad X\longmapsto \imath_X\omega_x
\end{equation}
from $\omega$ is injective, where $\imath_X\omega_x$ is the interior product of $\omega_x$ by $X$. 
\end{dfn}
A manifold equipped with an $E$-valued $n$-plectic form is called an $E$-valued $n$-plectic manifold. 
If a ${\rm d}^E_{\nabla}$-closed $E$-valued $(n+1)$-form $\omega$ is not necessarily non-degenerate, then it is called an $E$-valued pre-$n$-plectic form (or structure). 
A manifold equipped with an $E$-valued pre-$n$-plectic form is called an $E$-valued pre-$n$-plectic manifold. 
We denote by $(M,\omega,E,\nabla^E)$ an $E$-valued (pre-)$n$-plectic manifold $M$ together with an $E$-valued (pre-)$n$-plectic form $\omega\in \Omega^{n+1}(M,E)$ and a connection $\nabla^E$ on $E$. 

\begin{ex}[Pre-symplectic manifolds]\label{sec2:ex_symplectic}
If $M$ is a {\rm (}pre-{\rm )}symplectic manifold together with a {\rm (}pre-{\rm )}symplectic form $\omega$, we consider the trivial line bundle 
$\underline{\mathbb{R}}_M:=M\times \mathbb{R}$ with the trivial connection. 
Then, $\Omega^k(M,\underline{\mathbb{R}}_M)=\Omega^k(M)$ for $k\geqq 0$ and $\nabla\varphi={\rm d}\varphi$ for $\varphi\in\Omega^k(M)$. 
Therefore, $\omega$ is a $\underline{\mathbb{R}}_M$-valued {\rm (}pre-{\rm )}$1$-plectic form. 
\end{ex}

\begin{ex}\label{sec2:ex_family of symplectic}
Given a family of {\rm (}pre-{\rm )}symplectic structures $\{\omega_i\}_{i=1}^d$ on $M$, we define an $\mathbb{R}^d$-valued $2$-form $\omega^{(d)}$ on $M$ 
by $\omega^{(d)}:=\sum_{i=1}^d\omega_i\otimes \boldsymbol{e}_i$, where $\{\boldsymbol{e}_i\}_{i=1}^d$ denotes the standard basis of $\mathbb{R}^d$. 
$\omega^d$ is an $\underline{\mathbb{R}^d}_M$-valued {\rm (}pre-{\rm )}$1$-plectic form on $M$ with respect to the trivial connection 
$\mathbf{d}\eta = \sum_{i}({\rm d}\eta_i)\,\boldsymbol{e}_i$, where $\eta\in \Omega^k(M,\mathbb{R}^d)~(k\geqq 0)$. 
\end{ex}

\begin{ex}[Polysymplectic manifolds]
Let $k\geqq 1$ and $\{\boldsymbol{e}_i\}_{i=1}^k$ the standard basis of $\mathbb{R}^k$. 
A manifold $M$ equipped with a closed nondegenerate $\mathbb{R}^k$-valued $2$-form $\omega=\sum_i^k\omega_i\otimes \boldsymbol{e}_i$ 
is called a $k$-polysymplectic manifold. 
Every $k$-polysymplectic manifold $M$ is an $\underline{\mathbb{R}^k}_M$-valued $1$-plectic manifold with respect to 
the trivial connection $\mathbf{d}$. 
\end{ex}

\begin{ex}[Lie groups]\label{sec2:ex_Lie groups}
Let $G$ be a compact connected Lie group with a Lie algebra $\mathfrak{g}$, and $\lambda_L\,({\rm resp}. \lambda_R)$ the left {\rm (}{\rm resp}. right{\rm )}Maurer-Cartan form on $G$. 
That is, both $\lambda_L$ and $\lambda_R$ are $\mathfrak{g}$-valued 2-forms defined by 
\[
 (\lambda_L)_g(v_g) := (\mathrm{d}L_{g^{-1}})_g(v_g),\quad (\lambda_R)_g(v_g) := (\mathrm{d}R_{g^{-1}})_g(v_g)~;\quad g\in G, 
\]
where $L_g$ and $R_g$ denote the left translation and the right translation, respectively. 
Define a connection $\nabla^{\mathfrak{g}}$ 
on $\underline{\mathfrak{g}}_G:=G\times \mathfrak{g}$ by the trivial connection $\nabla^{\mathfrak{g}}f:={\rm d}f$, 
where $f$ is a section of $\underline{\mathfrak{g}}_G$ considered as a map from $G$ to $\mathfrak{g}$. 
Then, $G$ is a $\underline{\mathfrak{g}}_G$-valued pre-$1$-plectic manifold by ${\rm d}\lambda_L$ or ${\rm d}\lambda_R$ with respect to $\nabla^{\mathfrak{g}}$. 
\end{ex}

\begin{ex}[The curvature 2-form]
Let $M$ be a manifold, and $E$ a vector bundle equipped with a connection $\nabla^E$. 
Let $R^E_{\nabla}\in \Omega^2(M,{\rm End}\,E)$ be the curvature 2-form of $\nabla^E$. 
The connection $\nabla^E$ induces a connection on ${\rm End}\,E$ by 
\[
(\nabla^{{\rm End}E}_X\Phi)(e) := \nabla^E_X(\Phi(e)) - \Phi\bigl(\nabla^E_Xe\bigr), 
\]
where $\Phi\in {\rm End}\,E,\,X\in \mathfrak{X}(M)$ and $e\in \varGamma(E)$. 
Then, the 2-form $R^E_{\nabla}$ satisfies the identity ${\rm d}_{\nabla}^{{\rm End} E}R_{\nabla}^E=0$~{\rm (}the Bianchi identity{\rm )}. Therefore, $R_{\nabla}^E$ is an 
${\rm End} E$-valued pre-1-plectic form on $M$. 
\end{ex}

\begin{ex}
Let $M$ be a manifold, and $E$ a vector bundle equipped with a flat connection $\nabla$. Consider a vector bundle ${\rm Hom}(TM,E):=\coprod_{x\in M}{\rm Hom}(T_xM,E_x)\overset{\pi}{\to}M$, 
and define a 1-form $\vartheta$ on ${\rm Hom}(TM,E)$ by 
\[
\vartheta_{(x,\Phi)}(X):=\left(\Phi\circ({\rm d}\pi)_{(x,\Phi)}\right)(X), 
\]
where $X$ is any tangent vector to ${\rm Hom}(TM,E)$ at $(x,\Phi)\in {\rm Hom}(TM,E)$. 
Since $\nabla$ is a flat connection, an $E$-valued 2-form ${\rm d}^E_{\nabla}\vartheta\in \Omega^2({\rm Hom}(TM,E),E)$ is closed with respect to ${\rm d}^E_{\nabla}$. 
\end{ex}

\begin{ex}[Pre-$n$-plectic manifolds]
Similarly to the case of Example \ref{sec2:ex_symplectic}, a {\rm (}pre-{\rm )}$n$-plectic manifold $(M,\omega)$ is a $\underline{\mathbb{R}}_M$-valued {\rm (}pre-{\rm )}$n$-plectic manifold. 
\end{ex}
\bigskip 

Incidentally, we let $(M,\omega)$ be a (pre-)$n$-plectic manifold and assume that $M$ admits a Riemann metric $g$. 
The (pre-)$n$-plectic form $\omega$ can be thought of as a $T^*M$-valued $n$-form 
\[
\tilde{\omega}_x:T_xM\times \cdots\times T_xM \longrightarrow T_x^*M,\quad (v_1,\cdots,v_n)\longmapsto \imath_{v_1\wedge\cdots \wedge v_n}\omega_x, 
\]
where $x\in M$. The Riemannian metric $g$ induces the Levi-Civita connection $\nabla^g$ on $T^*M$ by 
\[
\left\langle\nabla_X^g\alpha, Y \right\rangle:= (\nabla_X^g\alpha)(Y) := X\bigl(\alpha(Y)\bigr) - \alpha\bigl(\nabla_X^gY\bigr) 
\]
for any $\alpha\in \Omega^1(M)$ and $X,Y\in \mathfrak{X}(M)$. 

\begin{prop}
Let $\omega$ be a {\rm (}pre-{\rm )}$n$-plectic form. The induced form $\tilde{\omega}$ is a 
$T^*M$-valued {\rm (}pre-{\rm )} $(n-1)$-plectic form with respect to the Levi-Civita connection $\nabla^g$ if and only if $\nabla^g\omega=0$. 
\end{prop}

\begin{proof}
For any vector fields $X_1,\cdots,X_{n+1},Y$ on $M$, a covariant exterior derivative ${\rm d}^g_\nabla$ of $\nabla^g$ is given by 
\begin{align*}
\bigl\langle({\rm d}^g_{\nabla}\tilde{\omega})(X_1,\cdots,X_{n+1}),\,Y\bigr\rangle 
&=  \sum_{i=1}(-1)^{i+1}\left\langle \nabla_{X_i}^g\bigl(\tilde{\omega}(X_1,\cdots,\check{X_i},\cdots,X_{n+1})\bigr),\,Y\right\rangle \\
&\quad+ \sum_{i<j}(-1)^{i+j}\left\langle \tilde{\omega}([X_i,X_j],X_1,\cdots,\check{X_i},\cdots,\check{X_j},\cdots,X_{n+1}),\,Y\right\rangle . 
\end{align*}
From the fact that the torsion of $\nabla^g$ is zero, it follows that 
\begin{align*}
\left\langle \nabla_{X_i}^g\bigl(\tilde{\omega}(X_1,\cdots,\check{X_i},\cdots,X_{n+1})\bigr),\,Y\right\rangle 
&= X_i\left(\omega(X_1,\cdots,\check{X}_i,\cdots,X_{n+1},Y)\right)\\
&\quad - \omega(X_1,\cdots,\check{X}_i,\cdots,X_{n+1},\nabla^g_YX_i)\\ 
&\qquad -(-1)^{n}\omega([X_i,Y],X_1,\cdots,\check{X}_i,\cdots,X_{n+1}) .
\end{align*}
Thus, ${\rm d}^g_{\nabla}\tilde{\omega}$ is calculated as 
\begin{align*}
&\bigl\langle({\rm d}^g_{\nabla}\tilde{\omega})(X_1,\cdots,X_{n+1}),\,Y\bigr\rangle \\
=\, & ({\rm d}\omega)(X_1,\cdots,X_{n+1},Y)\\
&\qquad + (-1)^nY\left(\omega(X_1,\cdots,X_{n+1})\right) 
+ \sum_{i=1}^{n+1}(-1)^i\omega(X_1,\cdots,\check{X}_i,\cdots,X_{n+1},\nabla^g_YX_i) \\
=\, & ({\rm d}\omega)(X_1,\cdots,X_{n+1},Y)\\
&\qquad + (-1)^n\left\{Y\left(\omega(X_1,\cdots,X_{n+1})\right) 
- \sum_{i=1}^{n+1}\omega(X_1,\cdots,\nabla^g_YX_i,\cdots,X_{n+1})\right\}\\
=\, & ({\rm d}\omega)(X_1,\cdots,X_{n+1},Y) + (-1)^n(\nabla^g_Y\omega)(X_1,\cdots,X_{n+1}). 
\end{align*}
Therefore, we have ${\rm d}^g_{\nabla}\tilde{\omega}={\rm d}\omega + (-1)^n\nabla^g\omega$. The assertion follows from this equation. 
\end{proof}
\bigskip 

In symplectic geometry, a vector field $X_f$ is associated with a function $f$ by ${\rm d}f=\imath_{X_f}\omega$. Such a function is said to be Hamiltonian, and 
$X_f$ is called the Hamiltonian vector field of $f$. It is checked easily that $X_f$ is a symplectic vector field, that is, $\mathcal{L}_{X_f}\omega = 0$. 
Similar notions are also defined in multi-symplectic geometry, being expanded to differential forms. Namely, if $\omega$ is a (pre-)$n$-plectic form and $\alpha$
is a $(n-1)$-form which satisfies ${\rm d}\alpha = \imath_{X_\alpha}\omega$ for some vector field $X_\alpha$, we say $\alpha$ is Hamiltonian and $X_\alpha$ is the 
Hamiltonian vector field corresponding to $\alpha$. 
A vector field $X$ is said to be multi-symplectic if $\mathcal{L}_X\omega=0$. By the Cartan formula, the Hamiltonian vector field $X_\alpha$ is also multi-symplectic. 

Taking those into account, we introduce an extended version of the Hamiltonian vector field into our discussion. 
We let $(M,\omega,E,\nabla^E)$ be an $E$-valued (pre-)$n$-plectic manifold. Given a vector field $X$ and an $E$-valued $k$-form $\varphi$ on $M$, we 
define the covariant Lie derivative of $\varphi$ with respect to $X$ by 
\begin{equation}\label{sec2:covariant Lie derivative}
(\mathcal{L}^{\nabla}_X\varphi)(X_1,\cdots,X_k) := \nabla^E_X\left(\varphi(X_1,\cdots,X_k)\right) - \sum_{i=1}^{k}\varphi\left(X_1,\cdots,[X,X_i],\cdots,X_k\right), 
\end{equation}
for $X_1,\cdots,X_k\in \mathfrak{X}(M)$. When $E=\underline{\mathbb{R}}_M$ and $\nabla^E ={\rm d}$, $\mathcal{L}_X^{\nabla}$ is none other than the Lie derivative of differential forms. 
\medskip 

The following lemma is proven in the same way as the proof of the Cartan formula.

\begin{prop}[Extended Cartan's formula]\label{sec2:extended Cartan}
Let $X,Y\in \mathfrak{X}(M)$ and $\varphi\in \Omega^k(M,E)$. Let $R^E_{\nabla}:={\rm d}^E_{\nabla}\circ \nabla^E\in \Omega^2(M,{\rm End}\,E)$ be the curvature of $\nabla^E$. 
Then, 
\begin{enumerate}[\quad \rm (1)]
\item $\mathcal{L}_X^{\nabla} = \imath_X\circ{\rm d}^E_{\nabla} + {\rm d}^E_{\nabla}\circ \imath_X.$
\item $\imath_{[X,Y]} = \mathcal{L}^{\nabla}_X\circ\imath_Y - \imath_Y\circ\mathcal{L}^{\nabla}_X$.
\item $\left(\mathcal{L}_X^{\nabla}\circ{\rm d}^E_{\nabla}\right)\varphi = \bigl(\imath_X R^E_{\nabla}\bigr)\wedge \varphi + \left({\rm d}^E_{\nabla}\circ \mathcal{L}^{\nabla}_X\right)\varphi$. 
Alternatively, 
\begin{align*}
\left(\mathcal{L}_X^{\nabla}\circ{\rm d}^E_{\nabla}\varphi\right)(X_1,\cdots, X_{k+1})&= 
\sum_{i=1}^{k+1}(-1)^{i+1}R^E_{\nabla}(X,X_i)\left(\varphi(X_1,\cdots,\check{X}_i,\cdots, X_{k+1})\right)\\
&\quad + \left({\rm d}^E_{\nabla}\circ\mathcal{L}^{\nabla}_X\varphi\right)(X_1,\cdots,X_{k+1}). 
\end{align*}
\end{enumerate}
\end{prop}

From Proposition \ref{sec2:extended Cartan} (3), one finds that the Lie derivative commute with the differential ${\rm d}$, 
while the covariant Lie derivative \eqref{sec2:covariant Lie derivative} does not do with ${\rm d}^E_{\nabla}$. 
It holds that $\mathcal{L}_X^{\nabla}\circ{\rm d}^E_{\nabla}={\rm d}^E_{\nabla}\circ \mathcal{L}_X^{\nabla}$ if $\nabla^E$ is a flat connection. 

\begin{dfn}
Let $(M,\omega,E,\nabla^E)$ be an $E$-valued {\rm (}pre-{\rm )}$n$-plectic manifold. An $E$-valued $(n-1)$-form $\varphi\in \Omega^{n-1}(M,E)$ is said to be a pseudo-Hamiltonian 
if there exists a vector field $X_{\varphi}\in \mathfrak{X}(M)$ such that ${\rm d}^E_{\nabla}\varphi = \imath_{X_{\varphi}}\omega$. 
We call $X_{\varphi}$ the pseudo-Hamiltonian vector field corresponding to $\varphi$. 
\end{dfn}

We denote by ${\rm pHam}^{n-1}(M,E,\nabla^E)$ the space of all pseudo-Hamiltonian ${\rm (}n-1{\rm )}$-forms on $M$. One might expect that $\mathcal{L}_{X_\varphi}^{\nabla}\omega$ 
vanishes for any $\varphi\in {\rm pHam}^{n-1}(M,E,\nabla^E)$, but this is in general not the case. 
Namely, even if $\varphi$ is a pseudo-Hamiltonian, $\mathcal{L}_{X_\varphi}^{\nabla}\omega=0$ does not necessarily hold 
because $\mathcal{L}_{X_\varphi}^{\nabla}\omega=\varphi\wedge R^E_{\nabla}$ from Proposition \ref{sec2:extended Cartan} (1). 
If $({\rm d}^E_{\nabla}\circ{\rm d}^E_{\nabla})\varphi=0$ is satisfied, then $\mathcal{L}_{X_\varphi}^{\nabla}\omega=0$ holds. 
\medskip 

\begin{dfn}\label{sec2:bracket}
For $\varphi,\psi\in {\rm pHam}^{n-1}(M,E,\nabla^E)$, we define an $E$-valued $(n-1)$-form $\{\varphi,\psi\}$ as 
\begin{equation*}
\{\varphi,\,\psi\}:= \imath_{X_{\psi}}\imath_{X_{\varphi}}\omega=\imath_{X_{\varphi}\wedge X_{\psi}}\omega \in \Omega^{n-1}(M,E). 
\end{equation*}
\end{dfn}
\medskip 

The operation is skew-symmetric and well-defined. 
Actually, if both $X_{\varphi}$ and $Y_{\varphi}$ are the pseudo-Hamiltonian vector fields for $\varphi\in {\rm pHam}^{n-1}(M,E,\nabla^E)$, then 
$\imath_{X_\psi}\imath_{X_\varphi}\omega = \imath_{X_\psi}{\rm d}^E_{\nabla}\varphi=\imath_{X_\psi}\imath_{Y_\varphi}\omega$. 
This means that the operation in Definition \ref{sec2:bracket} is well-defined.

\begin{lem}\label{sec2:E-hamiltonian}
Let $\varphi, \psi\in {\rm pHam}^{n-1}(M,E,\nabla^E)$. Then, 
\[
\imath_{[X_{\psi},X_{\varphi}]}\omega = {\rm d}^E_{\nabla}\{\varphi,\psi\} - \imath_{X_\varphi}(\psi\wedge R^E_{\nabla}) + \imath_{X_\psi}(\varphi\wedge R^E_{\nabla}).
\]
\end{lem}
\begin{proof}
The assertion follows from Proposition \ref{sec2:extended Cartan}, immediately. 
\end{proof}

By Lemma \ref{sec2:E-hamiltonian}, we find that the bracket $\{\varphi,\psi\}$ is {\em not} necessarily a pseudo-Hamiltonian form, in general. 
If $\nabla^E$ is a flat connection, $\{\varphi,\psi\}$ is pseudo-Hamiltonian again because the vector field $[X_{\psi},X_{\varphi}]$ satisfies 
$\imath_{[X_{\psi},X_{\varphi}]}\omega = {\rm d}^E_{\nabla}\{\varphi,\psi\}$. 

%
\section{$E$-valued Homotopy Momentum Sections}
\subsection{Differential Geometry of Lie Algebroids}

Let $M$ be a smooth manifold. 
A Lie algebroid over a smooth manifold $M$ is a smooth vector bundle $A\to M$ endowed with a bundle map $\rho:A\to TM$, called the anchor map, and the Lie bracket 
$[\cdot,\,\cdot]$ on each fiber which satisfy 
\[
[\alpha,\,f\beta] = f[\alpha,\,\beta] + \bigl(\rho(\alpha)f\bigr)\beta
\]
for any smooth section $\alpha,\,\beta \in \varGamma(A)$ and smooth function $f\in C^\infty(M)$. 
A Lie algebroid $A\to M$ is sometimes denoted by the tetrad $(A,\,M,\,\rho,\,[\cdot,\,\cdot])$. 
A tangent bundle of any manifold $M$ is a Lie algebroid whose the Lie bracket is the one for vector fields and the anchor is the identity map $\rho={\rm id}_{TM}$ on $TM$. 
This is called the tangent algebroid of $M$. 
For further examples of Lie algebroids, refer to \cite{DZpoi05, Mgen05} and the reference therein. 

Let $(A,\,M,\,\rho,\,[\cdot,\,\cdot])$ be a Lie algebroid. 
We call a section of the exterior bundle $\wedge^kA^*$ of $A$ with $k\geqq 0$ 
an $A$-differential $k$-form, and denote by $\Omega^k_{A}(M)$ the space $\varGamma(\wedge^k A^*)$ of $A$-differential $k$-forms. 
$\Omega_{A}^{\bullet}(M)$ yields the differential complex 
\[
\cdots\overset{{\eth}^A}{\longrightarrow} \Omega_{A}^{k-1}(M) \overset{{\eth}^A}{\longrightarrow} \Omega_{A}^k(M)\overset{{\eth}^A}{\longrightarrow} 
\Omega^{k+1}_{A}(M)\overset{{\eth}^A}{\longrightarrow} \cdots\, , 
\]
with the differential operator 
\begin{align*}
({\eth}^A\theta)(\alpha_1,\cdots,\alpha_{k+1}):= \sum_{i=1}^{k+1}(-1)^{i-1}&\rho(\alpha_i)\bigl(\theta(\alpha_1,\cdots,\check{\alpha}_i,\cdots,\alpha_{k+1})\bigr)\\
&+\sum_{i<j}(-1)^{i+j}\theta([\alpha_i,\,\alpha_j],\alpha_1,\cdots,\check{\alpha}_i,\cdots,\check{\alpha}_j,\cdots,\alpha_{k+1}), 
\end{align*}
where $\theta\in \Omega_{A}^k(M)$ and $\alpha_1,\cdots,\alpha_{k+1}\in \varGamma(A)$. Needless to say, the differential complex $(\Omega_{A}^{\bullet}(M),\,{\eth}_A)$ 
is none other than the de Rham complex $(\Omega^{\bullet}(M),\,{\rm d})$ in the case of $A=TM$. 

For non-negative integers $m, k$ such that $m\geqq k$, we introduce an operator 
\begin{equation}\label{sec3:general interior product}
\imath^k_\rho:\Omega^m(M)\to \Omega^{m-k}(M)\otimes_{C^\infty(M)} \Omega^k_{A}(M)
\end{equation} 
as 
\begin{equation}\label{sec2:general inner product 2}
(\imath^k_\rho\eta)(\alpha_1,\cdots,\alpha_k):= \imath_{\rho(\alpha_1)\wedge\cdots\wedge\rho(\alpha_k)}\eta \in \Omega^{m-k}(M),
\end{equation}
where $\eta\in \Omega^m(M)$ and $\alpha_1,\cdots,\alpha_k\in \varGamma(A)$. 
Particularly, $\imath_{\rho}^0$ is the identity map on $\Omega^m(M)$. 
It is checked easily that $\imath^1_{\rho}({\rm d}f)={\eth}^A f$ for $f\in C^\infty(M)$. In fact, we have 

\begin{lem}\label{sec3:lem_inner product}
It holds that $({\eth}^A\circ\imath_{\rho}^m)\,\eta = (\imath_{\rho}^{m+1}\circ {\rm d})\,\eta$ for any $\eta\in \Omega^m(M)$. 
\end{lem}
\medskip

Let $E$ be a vector bundle over the base manifold $M$ of $A$. We denote by $\Omega^k_{A}(M, E)$ the space of $A$-differential $k$-forms with values in $E$, that is, 
\[
\Omega_{A}^k(M,E):=\Omega_{A}^k(M)\otimes_{C^\infty(M)}\varGamma(E)\cong \varGamma(\wedge^k A^*\otimes E). 
\]

Similarly to the case of $A=TM$, one can define a connection of a Lie algebroid $A$ on $E$. 
An $A$-connection (or an $A$-covariant derivative) $\mho^E$ on $E$ is defined to be an $\mathbb{R}$-linear mapping 
\[
 \mho^E:\varGamma(E)\longrightarrow \Omega^1_{A}(M,E),\quad s\longmapsto \mho^E s
\]
which satisfies 
\[
 \mho^E(fs) = f\,\mho^E s + \imath^1_{\rho}({\rm d}f) \otimes s 
\]
for any $s\in \varGamma(E)$ and $f\in C^\infty(M)$. If $\alpha\in \varGamma(A)$, 
we often use the notation $\mho^E_{\alpha}s$ for $(\mho^E s)(\alpha)\in \varGamma(E)$. 
Given an $A$-connection, one can get an $\mathbb{R}$-linear operator ${\eth}_E^{\mho}:\Omega_{A}^k(M,E)\to \Omega_{A}^{k+1}(M,E)$ with $k\geqq 0$ by 
\[
 {\eth}_{\mho}^E(\theta\otimes s):= {\eth}^A\theta \otimes s + (-1)^k\theta\wedge \mho^E s \,;\quad \theta\in \Omega_{A}^k(M),\, s\in \varGamma(E), 
\]
and consequently, a sequence 
\begin{equation}\label{sec2;sequence}
0\longrightarrow \varGamma(E)\overset{\mho^E}{\longrightarrow} \Omega_{A}^1(M,E)\overset{{\eth}_{\mho}^E}{\longrightarrow} \Omega_{A}^2(M,E)
\overset{{\eth}_{\mho}^E}{\longrightarrow} \cdots\, .
\end{equation}
To be exact, (\ref{sec2;sequence}) is {\it not} a differential complex, that is, ${\eth}_{\mho}^E\circ {\eth}_{\mho}^E= 0$ fails to hold in general. 
However, it has the following properties;

\begin{lem}\label{sec2;lemma_covariant derivative} 
The operator ${\eth}^{\mho}_E$ satisfies the followings{\rm :} 
\begin{enumerate}[\quad \rm (1)] 
 \item ${\eth}_{\mho}^E=\mho^E$ when $k=0$. 
 \item ${\eth}_{\mho}^E(\theta\wedge \tau)={\eth}_A\theta\wedge \tau + (-1)^k\theta\wedge {\eth}_{\mho}^E\tau$, 
where $\theta\in \Omega_{A}^k(M)$ and $\tau\in \Omega_{A}^\ell(M,E)$.
\end{enumerate}
\end{lem}
\noindent 
The sequence (\ref{sec2;sequence}) is a Lie algebroid counterpart of the one from a covariant exterior derivative in the connection theory. 
So, we call ${\eth}_{\mho}^E$ the $A$-covariant exterior derivative of $\mho^E$. 
\medskip

The notion of a curvature of $\mho^E$ also can be considered by the same manner as the case of $A=TM$. The curvature of $\mho^E$, denoted by ${\Re}_{\mho}^E$, is 
an ${\rm End}_{C^\infty(M)}(E)$-valued $A$-differential 2-form 
\[
{\Re}_{\mho}^E : \varGamma(A)\times \varGamma(A) \longrightarrow {\rm End}_{C^\infty(M)}(\varGamma(E))
\] 
defined as
\[
{\Re}_{\mho}^E(\alpha,\beta)\,s
:= \bigl(\mho^E_\alpha\circ \mho^E_\beta - \mho^E_\beta\circ \mho^E_\alpha - \mho^E_{[\alpha, \beta]}\,\bigr)s \in \varGamma(E).
\]
Remark that ${\Re}_{\mho}^E$ can be identified with the $C^\infty(M)$-linear map ${\eth}_{\mho}^E\circ \mho^E:\varGamma(E)\to \Omega_{A}^2(M,E)$.

\subsection{Definition}

Let $E$ be a vector bundle over a manifold $M$ equipped with a vector bundle connection $\nabla^E:\varGamma(E)\to \Omega^1(M,E)$. Let $\omega$ be an $E$-valued (pre-)$n$-plectic form on $M$. 
Suppose that $(A,\,M,\,\rho,\,[\cdot,\,\cdot])$ is a Lie algebroid 
over $(M,\,\omega)$ with a vector bundle connection $\nabla^A:\varGamma(A)\to \Omega^1(M,A)$. 
We denote by $\Omega^{p,q}(M,A)^E$ the space of all smooth sections of $\wedge^pT^*M\otimes \wedge^qA^*\otimes E$ for integers $p\geqq 0,\,q\geqq 1$. Namely, 
\[
\Omega^{p,q}(M,A)^E := \Omega^p(M)\otimes_{C^\infty(M)}\Omega_{A}^q(M,E) \cong \varGamma(\wedge^pT^*M\otimes {\rm Hom}(\wedge^qA, E))\, . 
\]
We extend $\nabla^A$ to a vector bundle connection on ${\rm Hom}(\wedge^qA, E)\cong \wedge^q A^*\otimes E$ (which we denote by $\nabla^{A^*\otimes E}$) by 
\begin{equation}\label{sec3;eqn_connection}
\bigl({\nabla_X^{A^*\otimes E}}(\theta\otimes s)\bigr)(\alpha_1,\cdots,\alpha_q):=\nabla^E_X\bigl(\theta(\alpha_1,\cdots,\alpha_q)s\bigr) 
  - \sum_{i=1}^q\theta(\alpha_1,\cdots,\nabla^A_X\alpha_i,\cdots,\alpha_q)s, 
\end{equation}
where $\theta\otimes s \in\Omega_A^q(M,E)=\Omega_A^q(M)\otimes \varGamma(E)$, $\alpha_1,\cdots,\alpha_q\in \varGamma(A)$ and $X \in \mathfrak{X}(M)$. 
Then, we get a covariant exterior derivative of $\nabla^{A^*\otimes E}$
\begin{equation*}
 \mathrm{d}_{\nabla}^{A^*\otimes E}:\Omega^{p,q}(M,A)^E\longrightarrow \Omega^{p+1,q}(M,A)^E 
\end{equation*}
which is given by 
\begin{align}\label{sec3;eqn_covariant exterior derivative}
\bigl(\mathrm{d}_{\nabla}^{A^*\otimes E}(\eta\otimes\phi)\bigr)(X_1,\cdots,X_{p+1})
:=& \sum_{i=1}^{p+1}(-1)^{i-1}{\nabla}^{A^*\otimes E}_{X_i}\bigl(\eta(X_1,\cdots,\check{X}_i,\cdots,X_{p+1})\,\phi\bigr) \notag \notag \\
&+ \sum_{i<j}(-1)^{i+j}\eta([X_i,\,X_j],X_1,\cdots,\check{X}_i,\cdots,\check{X}_j,\cdots,X_{p+1})\,\phi,
\end{align}
where $\eta\in \Omega^p(M),\,\phi\in \Omega_A^q(M,E)$ and $X_1,\cdots,X_{p+1}\in \mathfrak{X}(M)$. 
Alternatively, ${\rm d}_{\nabla}^{A^*\otimes E}$ is expressed as 
\[
\mathrm{d}_{\nabla}^{A^*\otimes E}(\eta\otimes\phi)={\rm d}\eta\otimes \phi + (-1)^p\eta\wedge \nabla^{A^*\otimes E}\phi. 
\] 
In particular, $\mathrm{d}_{\nabla}^{A^*\otimes E}=\nabla^{A^*\otimes E}$ when $p=0$. 
Note that $\mathrm{d}_{\nabla}^{A^*\otimes E}$ is an operator which only increases one in the degree $p$ and does {\it not} necessarily satisfy 
$\mathrm{d}_{\nabla}^{A^*\otimes E}\circ\mathrm{d}_{\nabla}^{A^*\otimes E}=0$. 
\medskip 

On the other hand, we define an $A$-connection $\mho$ on $TM$ as 
\begin{equation*}\label{sec2;eqn_A-connection}
 \mho^{TM}_\alpha X := \rho(\nabla_X^A\alpha) + [\rho(\alpha),\,X] 
\end{equation*}
where $\alpha\in \varGamma(A),\,X\in \mathfrak{X}(M)$. 
We extend it to an $A$-connection on $\wedge^pT^*M\otimes E$ (which we denote by the letter $\mho^{T^*M\otimes E}$) by 
\begin{equation*}\label{sec2;extended A-connection}
\bigl(\mho^{T^*M\otimes E}_{\alpha}(\eta\otimes s)\bigr)(X_1,\cdots,X_p):= 
\nabla^E_{\rho(\alpha)}\bigl(\eta(X_1,\cdots,X_p) s\bigr) - \sum_{i=1}^p\eta(X_1,\cdots,\mho^{TM}_{\alpha}X_i,\cdots,X_p) s, 
\end{equation*}
where $\eta\otimes s \in\Omega^p(M)\otimes \varGamma(E)$. When $p=0$, we set $\mho^{T^*M\otimes E}_{\alpha} = \nabla^E_{\rho(\alpha)}$ for any $\alpha\in \varGamma(A)$. 
Consequently, we get an $A$-covariant exterior derivative of $\mho^{T^*M\otimes E}$ 
\[
{\eth}_{\mho}^{T^*M\otimes E}:\Omega^{p,q}(M,A)^E\longrightarrow \Omega^{p,q+1}(M,A)^E
\]
by 
\begin{align}\label{sec3;eqn_covariant exterior derivative2}
\bigl({\eth}_{\mho}^{T^*M\otimes E}(\theta\otimes\varphi)\bigr)(\alpha_1,\cdots,\alpha_{q+1})
:=& \sum_{i=1}^{q+1}(-1)^{i-1}\mho^{T^*M\otimes E}_{\alpha_i}(\theta(\alpha_1,\cdots,\check{\alpha}_i,\cdots,\alpha_{q+1})\,\varphi) \notag \\
+& \sum_{i<j}(-1)^{i+j}\theta([\alpha_i,\,\alpha_j],\alpha_1,\cdots,\check{\alpha}_i,\cdots,\check{\alpha}_j,\cdots,\alpha_{q+1})\,\varphi, 
\end{align}
where $\theta\in\Omega_A^q(M), \varphi\in \Omega^p(M,E)$ and $\alpha_1,\cdots,\alpha_{q+1}\in \varGamma(A)$. ${\eth}_{\mho}^{T^*M\otimes E}$ increases one in the degree $q$. 
As with ${\rm d}_{\nabla}^{A^*\otimes E}$, the relation ${\eth}_{\mho}^{T^*M\otimes E}\circ {\eth}_{\mho}^{T^*M\otimes E}=0$ fails to hold in general. 
\bigskip 

We extend the operator $\imath^k_\rho$ in (\ref{sec3:general interior product}) to 
\[
\imath^k_\rho:\Omega^m(M,E)\to \Omega^{m-k}(M,E)\otimes_{C^\infty(M)} \Omega^k_{A}(M)
\] 
by 
\begin{equation}
\bigl(\imath^k_\rho(\eta\otimes s)\bigr)(\alpha_1,\cdots,\alpha_k):= \imath_{\rho(\alpha_1)\wedge\cdots\wedge\rho(\alpha_k)}\eta\otimes s \in \Omega^{m-k}(M,E),
\end{equation}
where $\eta\in \Omega^m(M),\, s\in \varGamma(E)$ and $\alpha_1,\cdots,\alpha_k\in \varGamma(A)$. The following result is generalization of Lemma \ref{sec3:lem_inner product}. 

\begin{lem}
It holds that $({\eth}^E_{\mho}\circ\imath_{\rho}^m)\,\varphi = (\imath_{\rho}^{m+1}\circ {\rm d}^E_{\nabla})\, \varphi$ for any $\varphi\in \Omega^m(M,E)$, 
where ${\rm d}^E_{\nabla}$ is a covariant exterior derivative of the vector bundle connection $\nabla^E$ by \eqref{sec2:eqn_cov ext der}. 
\end{lem}

\begin{dfn}[Bundle-valued homotopy momentum section]\label{sec1:dfn_HMS}
Let $(M,\omega,E,\nabla^E)$ be an $E$-valued {\rm (}pre-{\rm )} $n$-plectic manifold. 
Let $A$ be a Lie algebroid over $(M,\omega,E,\nabla^E)$ equipped with a vector bundle connection $\nabla^A$. 
An $E$-valued homotopy momentum section with respect to $(A, \nabla^A)$ is a formal sum $\mu = \sum_{k=0}^{n-1}\mu_k$ with 
$\mu_k\in \Omega^{k,n-k}(M,A)^E\cong \Omega^k(M)\otimes \Omega_A^{n-k}(M,E)$ which satisfies 
\[
\bigl(\mathrm{d}_{\nabla}^{A^*\otimes E} + {\eth}_{\mho}^{T^*M\otimes E}\bigr)\, \mu = \sum_{k=0}^n(-1)^{n-k}\imath_{\rho}^{n+1-k}\omega. 
\]
We call $\mu_k$ the $k$-th component of $\mu$. 
\end{dfn}
\medskip 

The following diagram may help the readers understand Definition \ref{sec1:dfn_HMS}. In the diagram, each $k$-th component of $\mu$ belongs to $\Omega^{k,n-k}(M,A)^E$, and 
$\iota_\rho^k\omega$ appears in the spaces $\Omega^{n+1-k,k}(M,A)^E$. 

\begin{equation*}
\begin{CD}
 \Omega_A^{n+1}(M,E) \\
 @A{{\eth}_{\mho}^{T^*M\otimes E}}AA  \\
 \Omega_A^n(M,E) @>{\nabla^{A^*\otimes E}}>> \Omega^{1,n}(M,A)^E  \\
 @A{{\eth}_{\mho}^{T^*M\otimes E}}AA @A{{\eth}_{\mho}^{T^*M\otimes E}}AA  \\
 \vdots @>{\nabla^{A^*\otimes E}}>> \vdots @>{{\rm d}_{\nabla}^{A^*\otimes E}}>> \Omega^{n+1-k,k}(M,A)^E \\
 @A{{\eth}_{\mho}^{T^*M\otimes E}}AA @A{{\eth}_{\mho}^{T^*M\otimes E}}AA @A{{\eth}_{\mho}^{T^*M\otimes E}}AA \\
 \Omega_A^2(M,E) @>{\nabla^{A^*\otimes E}}>> \Omega^{1,2}(M,A)^E @>{{\rm d}_{\nabla}^{A^*\otimes E}}>> \dots @>{{\rm d}_{\nabla}^{A^*\otimes E}} >> \Omega^{n-1,2}(M,A)^E \\
 @A{{\eth}_{\mho}^{T^*M\otimes E}}AA @A{{\eth}_{\mho}^{T^*M\otimes E}}AA @A{{\eth}_{\mho}^{T^*M\otimes E}}AA @A{{\eth}_{\mho}^{T^*M\otimes E}}AA \\
 \Omega_A^1(M,E) @>{\nabla^{A^*\otimes E}}>> \Omega^{1,1}(M,A)^E @>{{\rm d}_{\nabla}^{A^*\otimes E}}>> \cdots @>{{\rm d}_{\nabla}^{A^*\otimes E}} >> 
 \Omega^{n-1,1}(M,A)^E @>{{\rm d}_{\nabla}^{A^*\otimes E}}>> \Omega^{n,1}(M,A)^E
\end{CD}
\end{equation*}
\medskip 

\begin{ex}[Momentum maps for symplectic manifolds]\label{sec3:example_momentum maps for symplectic manifolds}
Let $G$ be a compact and connected Lie group with Lie algebra $\mathfrak{g}$, and $(M,\omega)$ a symplectic manifold. 
As mentioned in Example \ref{sec2:ex_symplectic}, $(M,\omega)$ is thought of as a $\underline{\mathbb{R}}_M$-valued $1$-plectic manifold. 
Assume that $(M,\omega)$ admits a left Hamiltonian $G$-action $\Phi_g:M\to M~(g\in G)$ with a $G$-equivariant momentum map $J:M\to \mathfrak{g}^*$. 
The infinitesimal generator of the $G$-action 
\[
\left.{\xi}_M\right|_{x} := \left.\frac{d}{dt}\right|_{t=0}\!\!\Phi_{\exp(t\xi)}(x)\, ;\quad \xi\in \mathfrak{g},\, x\in M
\]
is a left $\mathfrak{g}$-action on $M$ and defines a Lie algebroid structure on the trivial bundle $\mathfrak{g}\times M$ over $M$ by 
\[
[\alpha,\beta](x) := [\alpha(x),\,\beta(x)] - \alpha(x)_M|_x\beta + \beta(x)_M|_x\alpha 
\] 
for any section $\alpha,\beta\in \varGamma(\mathfrak{g}\times M)$. Here, we identify $\alpha,\,\beta$ with functions from $M$ to $\mathfrak{g}$. 
This algebroid is called the action Lie algebroid and denoted by $\mathfrak{g}\ltimes M$. 
Define a connection on the algebroid $A=\mathfrak{g}\ltimes M$ as the trivial connection $\nabla \alpha:=\mathbf{d}\alpha$. 
Regarding $J$ as a section of the dual bundle $\mathfrak{g}^*\times M$ and confining sections of $A$ to the constant ones, 
we find that the condition that ${\rm d}_{\nabla}J=\iota_{\rho}^1\omega$ is equivalent to ${\rm d}J^{\xi}=\imath_{\xi_M}\omega$, where $J^{\xi}$ is a function by 
$J^{\xi}(x):=\langle J(x),\xi\rangle$. 
And moreover, $J$ is equivariant if and only if ${\eth}_{\mho}J=-\iota_{\rho}^2\omega$. So, the equivariant momentum map for Hamiltonian action is thought of as a 
$\underline{\mathbb{R}}_M$-valued homotopy momentum section with respect to $(\mathfrak{g}\ltimes M, \mathbf{d})$. 
\end{ex}

\begin{ex}[Momentum sections]\label{sec3:example_momentum sections}
Let $M$ be a {\rm (}pre-{\rm )}symplectic manifold and $A$ a Lie algebroid over $M$ with a vector bundle connection $\nabla^A$. 
A bracket-compatible $\nabla^A$-momentum section is an $E$-valued homotopy momentum section for the case where $E$ is the trivial line bundle $E=\underline{\mathbb{R}}_M$ 
with a trivial connection $\nabla^E={\rm d}$. 
Here, we do not assume that the Lie algebroid $A$ is not presymplectically anchored with respect to $\nabla^A$, that is, the condition $R^A_{\nabla}\mu =0$ is {\em not} necessarily satisfied. 
\end{ex}

\begin{ex}[Homotopy momentum sections]
Let $M$ be a {\rm (}pre-{\rm )}$n$-plectic manifold and $A$ a Lie algebroid over $M$ with a vector bundle connection $\nabla^A$. 
A homotopy momentum section with respect to $\nabla^A$ is an $\underline{\mathbb{R}}_M$-valued homotopy momentum section. 
Remark that the sign convention in the right-hand side in the definition differs from the original one in \cite{HIhom22}.   
\end{ex}

\begin{ex}[Adjoint operators]
A Lie group $G$ is a $\mathfrak{g}$-valued $1$-plectic manifold by the derivative $\mathrm{d}\lambda_R$ of the right Maurer-Cartan form 
$\lambda_R$~{\rm (}see Example \ref{sec2:ex_Lie groups}{\rm )}. 
Consider the action algebroid $A=G\rtimes \mathfrak{g}$ by the right Lie algebra action $\mathfrak{g}\ni \alpha\mapsto \alpha^L\in \mathfrak{X}(G)$, where $\alpha^L$ denotes the left invariant 
vector field on $G$. 
For $\alpha\in \mathfrak{g}$, we define a $\mathfrak{g}$-valued map $\mu^{\alpha}$ on $G$ by $\mu^{\alpha}(h):=-\mathrm{Ad}_g\alpha$. 
By using the fundamental formula $\mathbf{d}\lambda_R-[\lambda_R,\,\lambda_R]_{\mathfrak{g}}=0$ and the fact that the adjoint representation $\mathrm{Ad}_h$ preserves the Lie bracket 
$[\cdot,\,\cdot]_{\mathfrak{g}}$ on $\mathfrak{g}$, we have 
\[
(\mathbf{d}\mu^{\alpha})_h(X_h) = (\mathbf{d}\lambda_R)_h(\alpha^L_h,\,X_h)\quad \text{and}\quad 
(\mathbf{d}\mu)_h(\alpha,\beta) = -(\mathbf{d}\lambda_R)_h(\alpha^L_h,\, \beta^L_h), 
\]
where $h\in G$ and $\alpha,\,\beta \in C^\infty(G,\mathfrak{g})$, considered as the constant functions. This means that $\mu$ is a $\mathfrak{g}$-valued homotopy momentum section on $G$. 
\end{ex}
\section{Quaternionic K\"{a}hler Symmetries}

\subsection{Quaternionic K\"{a}hler manifolds}

We begin with the section by recalling the fundamentals of a quaternionic K\"{a}hler manifold. 
Let $(M,g)$ be a Riemannian manifold of dimension $4m~(m\geqq 1)$. 
Suppose that there is a subbundle $\mathcal{Q}\subset {\rm End}_{C^\infty(M)}TM$ satisfying the following condition: 
at each point $x\in M$, there are a local coordinate neighborhood $U$ of $x$ and a local frame 
$\{J_1,J_2,J_3\}$ of $\mathcal{Q}|_U$ satisfying 
\begin{equation}\label{sec4:eqn_quaternionic relations}
J_aJ_b=-\delta_{ab}~{\rm id} + \varepsilon_{abc}J_c~;~\quad a,b,c=1,2,3. 
\end{equation}
Then, $(M,g,\mathcal{Q})$ is called an almost quaternionic manifold. 
Assume that the metric $g$ satisfies that $g(JX,JY)=g(X,Y)$ for any $J\in \varGamma(\mathcal{Q})$ and $X,Y\in \mathfrak{X}(M)$. 
Then, $\mathcal{Q}$ is embedded isometrically into $\wedge^2T^*M$ by $J\mapsto \omega_J$, where $\omega_J$ is a non-degenerate 2-form by $\omega_J(X,Y):=g(JX,Y)$. 
Consequently, one gets has a local frame $\{\omega_a\}$ of $\mathcal{Q}$ given by $\omega_a(X,Y):=g(J_aX,Y)$ and 
a global 4-form $\Theta^{\wedge}$, called the fundamental $4$-form 
\begin{equation*}
\Theta^{\wedge} = \sum_{i=1}^3\omega_i\wedge \omega_i. 
\end{equation*}
Since each $\omega_i$ is non-degenerate, $\Theta^{\wedge}$ is also non-degenerate. 
\begin{dfn}
An almost quaternionic manifold $(M,g,\mathcal{Q})$ of dimension greater than $4$ is called a quaternionic K\"{a}hler manifold if $\nabla^g\Theta^{\wedge}=0$, where 
$\nabla^g$ is the Levi-Civita connection. 
\end{dfn}
The condition $\nabla^g\Theta^{\wedge}=0$ in the definition implies that ${\rm d}\Theta^{\wedge}=0$ and each $\nabla^g\omega_i~(i=1,2,3)$ satisfies 
\begin{equation*}
\nabla^g_X\omega_j = \sum_{i=1}^3\alpha_{ij}(X)\,\omega_i, \quad X\in\mathfrak{X}(M),
\end{equation*}
where $\alpha_{ij}$ are 1-forms with $\alpha_{ij}=-\alpha_{ji}$. 
That is, the Levi-Civita connection $\nabla^g$ on $M$ preserves the subbundle $\mathcal{Q}$. 
If ${\rm dim}\,M=4$, then $\nabla^g\Theta^{\wedge}={\rm d}\Theta^{\wedge}=0$ is automatically satisfied. 
As for the case where $M$ is of four-dimension, a quaternionic K\"{a}hler manifold is defined to be an oriented Riemannian manifold which is both Einstein and self-dual. 
Any quaternionic K\"{a}hler manifold is a 3-plectic manifold by the fundamental $4$-form. 
For further details, refer to \cite{Iqua74,BGsas08} for example and the references therein. 

\begin{lem}\label{sec4:lem_quaternionic}
A quaternionic K\"{a}hler manifold is a $\mathcal{Q}$-valued $1$-plectic manifold. 
\end{lem}
\begin{proof}
Let $(M,g,\mathcal{Q})$ be a quaternionic K\"{a}hler manifold. We define a covariant tensor field $\Theta$ as 
\[
\Theta = \sum_{i=1}^3\omega_i\otimes \omega_i,
\]
where each $\omega_i$ is a local frame of the subbundle $\mathcal{Q}$. $\Theta$ is thought of as a $\mathcal{Q}$-valued $2$-form on $M$. 
It suffices to show that $\Theta$ is closed under the covariant exterior differential ${\rm d}^g_{\nabla}$ of $\nabla^g$. By applying ${\rm d}^g_{\nabla}$ to $\Theta$, it is calculated as
\begin{align*}
{\rm d}^g_{\nabla}\Theta &= \sum_{j=1}^3{\rm d}_g^{\nabla}(\omega_j\otimes \omega_j) = \sum_{j=1}^3({\rm d}\omega_j\otimes \omega_j + \omega_j\wedge \nabla^g\omega_j)\\
&= \sum_{j=1}^3\biggl\{{\rm d}\omega_j\otimes \omega_j + \sum_{i=1}^3(\omega_j\wedge\alpha_{ij})\otimes \omega_i\biggr\}\\
&= \sum_{j=i}^3\biggl({\rm d}\omega_i + \sum_{j=1}^3\alpha_{ij}\wedge \omega_j\biggr)\otimes \omega_i
\end{align*}
By the structure equations ${\rm d}\omega_i + \sum_{j=1}^3\alpha_{ij}\wedge\omega_j=0~(j=1,2,3)$, we have ${\rm d}^g_{\nabla}\Theta = 0$. 
It follows from the non-degeneracy of each $\omega_i$ that $\Theta$ is non-degenerate. Thus, the assertion is proved. 
\end{proof}

\subsection{Relation to a momentum map for a quaternionic K\"{a}hler manifold}

K. Galicki and H. B. Lawson devised a momentum mapping for quaternionic K\"{a}hler manifolds to discuss a reduction procedure for them in \cite{GLqua88}. 
They showed that if a quaternionic K\"{a}hler manifold $M$ has a non-zero scalar curvature, for each Killing vector field $V$ satisfying $\mathcal{L}_V\Theta^{\wedge}=0$, 
there exists a section $f_V\in \varGamma(\mathcal{Q})$ such that $\nabla f_V = \Theta_V$, where 
\[
\Theta_V := \sum_{i=1}^3(\imath_V\omega_i)\otimes \omega_i \in \Omega^1(M,\mathcal{Q}). 
\]
A Killing vector field $V$ with the condition that $\mathcal{L}_V\Theta^{\wedge}=0$ is called the quaternionic K\"{a}hler Killing vector field. 
The set of all quaternionic K\"{a}hler Killing vector fields is a Lie algebra by the natural Lie bracket for vector fields. 

\begin{dfn}[\cite{GLqua88}]
Let $\mathfrak{K}$ be the Lie subalgebra of a quaternionic K\"{a}hler Killing vector fields. 
The momentum map for $M$ is a section $f$ of $\mathfrak{K}^*\otimes \mathcal{Q}$ satisfying $\nabla^g f_V = \Theta_V$ for all $V\in \mathfrak{K}$. 
\end{dfn}
\medskip 

Here, the question arises: Is it possible to understand the momentum map $f$ of Galicki-Lawson's within the framework of the $E$-valued homotopy momentum section ? 

To answer it, we let $(M,g,\mathcal{Q})$ be a quaternionic K\"{a}hler manifold of dimension $4m~(m>1)$ with non-zero scalar curvature. 
Remark that the subbundle $\mathcal{Q}$ is equipped with the Levi-Civita connection $\nabla^g$. 
From Lemma \ref{sec4:lem_quaternionic}, $M$ is a $\mathcal{Q}$-valued 1-plectic manifold by $\Theta$. 
Assume that there is a vector subbundle $K\subset TM$ whose sections 
are quaternionic K\"{a}hler Killing vector fields.  Since $[V_1,V_2]\in \varGamma(K)$ for any $V_1,V_2\in \varGamma(K)$, 
$K$ is endowed with a Lie algebroid structure over $M$ whose anchor is the identity map ${\rm id}:K\to K\subset TM$. 

Define the $A$-connection $\mho^{TM}$ on $TM$ by 
\[
 \mho_V^{TM}X:= \nabla^g_XV + [V, X],\qquad V,\,X\in \mathfrak{X}(M). 
\]
However, since $\nabla^g$ is torsion-free, we see that $\mho^{TM}=\nabla^g$. 
Accordingly, the covariant exterior derivative ${\eth}_{\mho}^{T^*M\otimes \mathcal{Q}}: 
\Omega^{p,q}(M,K)^{\mathcal{Q}}\longrightarrow \Omega^{p,q+1}(M,K)^{\mathcal{Q}}$ of $\mho$ is given by 
\[
 {\eth}_{\mho}^{T^*M\otimes \mathcal{Q}}(\varphi\otimes\theta) = \varphi\otimes {\rm d}\theta + \nabla^g\varphi\wedge \theta,
\]
where $\varphi\in \Omega^p(M,\mathcal{Q}),\, \theta\in \Omega_{K}^q(M)$. 

\begin{thm}\label{sec4:thm_G-L}
A momentum map $f$ of Galicki-Lawson's is a $\mathcal{Q}$-valued homotopy momentum section with respect to the Levi-Civita connection $\nabla^g$ 
if and only if $f$ satisfies 
\begin{equation}\label{sec4:eqn_thm G-L}
f_{[V_1,V_2]} = - \sum_{i=1}^3\omega_i(V_1,V_2)\,\omega_i
\end{equation}
for any $V_1,V_2\in \varGamma(K)$. 
\end{thm}

\begin{proof}
Let $f$ be a momentum map of Galicki-Lawson's and $V\in\varGamma(K)$, the Lie subalgebra of a quaternionic K\"{a}hler Killing vector fields. 
Recall that $f$ is an element of $\Omega^1_{K}(M,\mathcal{Q})=\varGamma(K^*\otimes\mathcal{Q})$ and $f_V=f(V)\in \varGamma(Q)$. 
First, we shall show that the condition $\nabla^g f_V = \Theta_V$ is equivalent to $\nabla^g f=\imath_{\rm id}^1\Theta$. 
If $X\in \mathfrak{X}(M)$, then
\begin{equation}\label{sec4:eqn_thm}
(\nabla^g_X f)(V) = \nabla^g_Xf_V - f_{\nabla^g_XV}. 
\end{equation}
Since $\nabla^g$ is a metric connection and its torsion is vanishing, we have 
\[
(\mathcal{L}_Vg)(X,Y)=g(\nabla^g_XV,Y) + g(X,\nabla^g_YV)=0, 
\]
where $X,Y\in \mathfrak{X}(M)$. 
By setting $X=Y$, we get $\nabla^g_XV=0$. 
From this and \eqref{sec4:eqn_thm} it follows that $(\nabla^g f)(V)=\nabla^gf_V$. 
On the other hand, it is checked easily that $\Theta_V=(\imath^1_{\rm id}\Theta)(V)$. 
Therefore, $f$ is a momentum map of Galicki-Lawson's if and only if $f$ satisfies $\nabla^g f=\imath_{\rm id}^1\Theta$. 

Next, the covariant derivative of $f\in \Omega^1_K(M,\mathcal{Q})$ is calculated as 
\begin{align*}
({\eth}^{T^*M\otimes \mathcal{Q}}_{\mho}f)(V_1,V_2) &= \nabla^g_{V_1}f_{V_2} - \nabla^g_{V_2}f_{V_1} - f_{[V_1,V_2]}\\
&=(\nabla^g f_{V_2})(V_1) - (\nabla^g f_{V_1})(V_2) - f_{[V_1,V_2]}\\
&=\Theta_{V_2}(V_1) - \Theta_{V_1}(V_2) - f_{[V_1,V_2]}\\
&= -2\sum_{i=1}^3\omega_i(V_1,V_2)\,\omega_i - f_{[V_1,V_2]}.
\end{align*}
If the condition \eqref{sec4:eqn_thm G-L} is satisfied, then we have $({\eth}^{T^*M\otimes \mathcal{Q}}_{\mho}f)(V_1,V_2)=-\imath_{V_2}\imath_{V_1}\Theta$, which leads us 
to the condition ${\eth}^{T^*M\otimes \mathcal{Q}}_{\mho} f = -\imath_{\rm id}^2\Theta$. 
Conversely, if ${\eth}^{T^*M\otimes \mathcal{Q}}_{\mho} f = -\imath_{\rm id}^2\Theta$, \eqref{sec4:eqn_thm G-L} immediately follows from the above calculation. 
\end{proof}

The condition $\nabla^gf_V=\Theta_V$ indicates that a section $f_V$ of $\mathcal{Q}$ is a pseudo-Hamiltonian $0$-form corresponding to $V\in \mathfrak{K}$. 
Accordingly, we can define a bracket $\{f_{V_1},\,f_{V_2}\}$ for $f_{V_1}$ and $f_{V_2}$ by Definition \ref{sec2:bracket}. 
Then, the condition \eqref{sec4:eqn_thm G-L} is alternatively described as 
\begin{equation}\label{sec4:eqn2_thm}
f_{[V_1,V_2]} = - \{f_{V_1},\,f_{V_2}\}. 
\end{equation}
That is, the map $f\in {\rm Hom}(\mathfrak{K},\mathcal{Q})$ preserves the brackets up to the negative sign. In other words, 
$\{f_{V_1},\,f_{V_2}\}$ is a pseudo-Hamiltonian $0$-form again because 
\[
\nabla^{g}\{f_{V_1},f_{V_2}\}=\nabla^gf_{[V_1,V_2]}=\Theta_{[V_1,V_2]}=\imath_{[V_1,V_2]}\Theta. 
\] 
\subsection{Hyper-K\"{a}hler symmetries}

Momentum maps for hyper-K\"{a}hler manifolds have also been studied in \cite{HKLR87, Ggene87}, for instance. 
Recall that a hyper-K\"{a}hler manifold is a Riemannian manifold $(M, g)$ with three independent K\"{a}hler structures. 
More precisely, $(M, g)$ is a $4n$-dimensional Riemannian manifold endowed with three almost complex structures $I_1,\,I_2,\,I_3$ satisfying the same relation 
as \eqref{sec4:eqn_quaternionic relations} such that the $2$-forms $\omega_a(X,Y):=g(I_aX,Y)~(a=1,2,3)$ associated to both $I_a$ and $g$ are symplectic structures on $M$. 
As mentioned in Example \ref{sec2:ex_family of symplectic}, every hyper-K\"{a}hler manifold $(M,\,g)$ is considered as an $\mathbb{R}^3$-valued $1$-plectic manifold 
$(M,\,\omega^{(3)},\,\underline{\mathbb{R}^3}_M,\,\mathbf{d})$ . 

Let $G$ be a compact Lie subgroup of the isometry group acting properly and freely on $M$ from the left. 
Suppose that the $G$-action is Hamiltonian with respect to each of those symplectic forms $\omega_a$. Denote by $\mu_a$ each $G$-equivariant momentum map associated to $\omega_a$. 
Such a $G$-action is said to be hyper-hamiltonian \cite{BGsas08}. 
A hyper-K\"{a}hler momentum map is defined to be a map $\mu:M\to \mathfrak{g}^*\otimes\mathbb{R}^3$ by $\mu =\sum_{a}\mu_a\otimes\boldsymbol{e}_a$, 
where $\{\boldsymbol{e}_1,\boldsymbol{e}_2,\boldsymbol{e}_3\}$ is the standard basis of $\mathbb{R}^3$, and where $\mathfrak{g}^*$ denotes the dual of Lie algebra $\mathfrak{g}$ of $G$. 
\medskip 

Given a hyper-K\"{a}hler momentum map $\mu$ by a hyper-hamiltonian $G$-action, we consider the trivial bundle over $E=\underline{\mathbb{R}^3}_M$ and the action Lie algebroid 
$A=\mathfrak{g}\ltimes M$. 
We define a connection on each of them as the trivial one: $\nabla^{E}=\nabla^A=\mathbf{d}$. 
Each component $\mu_a$ of $\mu$ satisfies that 
\[
\imath_{\xi_M}\omega_a = {\rm d}(\mu_a^{\xi}) = (\mathrm{d}\mu_a)^{\xi} ;\quad a=1,2,3
\] 
for any $\xi\in \mathfrak{g}$~(see Example \ref{sec3:example_momentum maps for symplectic manifolds}). 
Here, $\xi$ is regarded as the constant section of $A$. 
This means that the same equation as 
$\mathrm{d}^{\mathfrak{g}^*\otimes\mathbb{R}^3}_{\nabla}\,\mu = \imath^1_{\rho}\omega$ in Definition \ref{sec1:dfn_HMS} holds for the hyper-K\"{a}hler momentum map $\mu$. 
Furthermore, from the condition that each $\mu_a$ is $G$-equivariant, we have 
\begin{equation*}
\mu^{[\alpha,\,\beta]} = - \{\mu^{\alpha},\,\mu^{\beta}\}  
\end{equation*}
for any $\alpha,\beta\in \mathfrak{g}$. 
This equation is equivalent to that $\eth_{\mho}^{T^*M\otimes \mathbb{R}^3}\mu = -\imath_{\rho}^2\omega^{(3)}$. 
Consequently, we get the following result by the same manner as the proof of Theorem \ref{sec4:thm_G-L}. 

\begin{prop}\label{sec4:prop_hyperkahler momentum maps}
Let $\mu:M\to \mathfrak{g}^*\otimes\mathbb{R}^3$ be a hyper-K\"{a}hler momentum map. When if we confine sections of the action Lie algebroid $\mathfrak{g}\ltimes M$ to 
the constant ones, $\mu$ is an $\mathbb{R}^3$-valued homotopy momentum section with respect to $(\mathfrak{g}\ltimes M,\, \mathbf{d})$. 
\end{prop}

\section{Compatibility with a Lie Algebroid}

For the subsequent discussion, we beforehand introduce the operation $\imath_\alpha$ for a section $\alpha\in \varGamma(A)$ by 
\[
\bigl(\imath_\alpha\nu\bigr)(X_1,\cdots,X_{n-1}):=\nu^{\alpha}(X_1,\cdots,X_{n-1})
:=\left\langle \nu(X_1,\cdots,X_{n-1}),\,\alpha\right\rangle\in \varGamma(E), 
\]
where $X_1,\cdots, X_{n-1} \in\mathfrak{X}(M)$ and $\nu\in \Omega^{n-1}(M)\otimes \Omega_A^{1}(M,E)$. 
We sometimes write $\imath_{\alpha}\nu$ as $\nu^\alpha$, for simplicity. 
\smallskip 

Let $n\geqq 1$. Let $(M,\omega,E,\nabla^E)$ be an $E$-valued (pre-)$n$-plectic manifold and $A$ a Lie algebroid over $M$ equipped with a vector bundle connection $\nabla^A$. 
Suppose that there exists an $E$-valued homotopy momentum section $\mu=\sum_{k=0}^{n-1}\mu_k$ with respect to $\nabla^{A}$. Recall that each component $\mu_k$ 
is an element in $\Omega^{k}(M)\otimes \Omega_A^{n-k}(M,E)$. 
\medskip 

As mentioned in the proof of Theorem \ref{sec4:thm_G-L}, a momentum map $f$ of Galicki-Lawson's satisfies the relation $(\nabla^g_X f)(V) = \nabla^g_Xf_V$. 
This means that the covariant derivative $\nabla^g$ and $\imath_V$ are commutative with each other. 
However, in general, the $(n-1)$-component $\mu_{n-1}$ of an $E$-valued homotopy momentum section does {\em not} necessarily satisfy that 
$(\imath_{\alpha}\circ{\rm d}^{\nabla})\mu_{n-1} = ({\rm d}^{\nabla}_E\circ\imath_{\alpha})\mu_{n-1}$ for any $\alpha\in\varGamma(A)$. 

\begin{dfn}\label{sec5:dfn_compatibility with A}
We say that the ${\rm (}n-1{\rm )}$-th component $\mu_{n-1}$ of $\mu$ is compatible with $A$ 
if it holds that $(\imath_{\alpha}\circ{\rm d}^{\nabla})\mu_{n-1} = ({\rm d}^{\nabla}_E\circ\imath_{\alpha})\mu_{n-1}$ for any $\alpha\in\varGamma(A)$.  
\end{dfn}

If the ${\rm (}n-1{\rm )}$-th component $\mu_{n-1}$ is compatible with $A$, it must simultaneously satisfy the condition 
${\rm d}^{\nabla}\mu_{n-1}=\imath_{\rho}^1\omega$ (see Definition \ref{sec1:dfn_HMS}). 
Then, by \eqref{sec3;eqn_connection} and \eqref{sec3;eqn_covariant exterior derivative}, we have 
\begin{align*}
\bigl((\imath_{\alpha}\circ{\rm d}^{\nabla})\mu_{n-1}\bigr)(X_1,\cdots,X_n)  
&=\sum_{i=1}^n(-1)^{i+1}\nabla^E_{X_i}\left({\mu_{n-1}}^{\alpha}(X_1,\cdots,\check{X}_i,\cdots,X_n)\right) \\
&\quad - \sum_{i=1}^n(-1)^{i+1}(\mu_{n-1})^{\nabla^A_{X_i}\alpha}(X_1,\cdots,\check{X}_i,\cdots,X_n)\\
&\qquad + \sum_{i<j}(-1)^{i+j}{\mu_{n-1}}^{\alpha}([X_i,X_j],X_1,\cdots,\check{X}_i,\cdots,\check{X}_j,\cdots,X_n). 
\end{align*}
for any $\alpha\in\varGamma(A)$. 

On the other hand, the covariant exterior derivative ${\rm d}^{\nabla}_E$ of $\imath_{\alpha}\mu_{n-1}\in \Omega^{n-1}(M,E)$ is calculated as
\begin{align*}
\left(({\rm d}^{\nabla}_E\circ \imath_{\alpha}){\mu_{n-1}}\right)(X_1,\cdots,X_n)&= \sum_{i=1}^n(-1)^{i+1}\nabla^E_{X_i}\left({\mu_{n-1}}^{\alpha}(X_1,\cdots,\check{X}_i,\cdots,X_n)\right) \\
&\quad + \sum_{i<j}(-1)^{i+j}{\mu_{n-1}}^{\alpha}([X_i,X_j],X_1,\cdots,\check{X}_i,\cdots,\check{X}_j,\cdots,X_n). 
\end{align*}
From those two formulas, we get directly the following proposition:

\begin{prop}\label{sec5:prop_momentum form}
The ${\rm (}n-1{\rm )}$-th component $\mu_{n-1}$ of a $E$-valued homotopy momentum section is compatible with $A$ if and only if it holds that
\[
\sum_{i=1}^n(-1)^{i+1}(\mu_{n-1})^{\nabla^A_{X_i}\alpha}(X_1,\cdots,\check{X}_i,\cdots,X_n) = 0 
\]
for any $\alpha\in\varGamma(A)$ and $X_1,\cdots, X_n \in\mathfrak{X}(M)$. 
\end{prop}

Obviously, the compatible component $\mu_{n-1}$ satisfies the relation ${\rm d}_E^{\nabla}({\mu_{n-1}}^{\alpha})=\imath^1_{\rho(\alpha)}\omega$ 
for each $\alpha\in \varGamma(A)$. 
This implies that ${\mu_{n-1}}^{\alpha}\in \Omega^{n-1}(M,E)$ is a pseudo-Hamiltonian $(n-1)$-form whose pseudo-Hamiltonian vector field is $\rho(\alpha)$. 
Namely, the compatibility is associated with the characteristic of being pseudo-Hamiltonian. 

\medskip 

Now, let us consider the case for $n=1$. Namely, $(M,\omega,E,\nabla^E)$ is an $E$-valued (pre-)$1$-plectic manifold and an $E$-valued homotopy momentum section 
$\mu$ is an element in $\varGamma(A^*\otimes E)$. 
Suppose that $\mu$ is compatible with $A$. Then, it satisfies that 
\[
 \forall\alpha\in\varGamma(A)\, ;\, (\nabla^{A^*\otimes E}\mu)^{\alpha} = {\nabla}^E{\mu}^{\alpha}=\imath_{\rho(\alpha)}\omega\quad \text{and}\quad {\eth}^{\mho}\mu = -\imath^2_{\rho}\omega. 
\]
$\mu^\alpha\in \varGamma(E)$ being a pseudo-Hamiltonian $0$-form, we can define a section of $E$ as 
\begin{equation}\label{sec5:eqn_bracket}
\{\mu^\alpha,\mu^\beta\} := \omega(\rho(\alpha),\,\rho(\beta))
\end{equation} 
for any $\alpha,\beta\in \varGamma(A)$~(see Definition \ref{sec2:bracket}).
From \eqref{sec3;eqn_covariant exterior derivative2} and $\mho_{\alpha}=\nabla^E_{\rho(\alpha)}$, it follows that
\begin{align*}
({\eth}^{\mho}\mu)(\alpha,\beta) 
= \nabla^E_{\rho(\alpha)}\mu^{\beta} - \nabla^E_{\rho(\beta)}\mu^{\alpha} - \mu([\alpha,\beta]).
\end{align*}
By using \eqref{sec5:eqn_bracket}, the first term $\nabla^E_{\rho(\alpha)}\mu^{\beta}$ in the right-hand side is calculated as
\[
\nabla^E_{\rho(\alpha)}\mu^{\beta} = (\imath_{\rho(\beta)}\omega)(\rho(\alpha))
=\omega(\rho(\beta),\rho(\alpha)) = -\{\mu^\alpha,\,\mu^\beta\}. 
\]
In a similar way, $\nabla^E_{\rho(\beta)}\mu^{\alpha}=\{\mu^\alpha,\,\mu^\beta\}$. 
From those relations, we obtain the following proposition. 

\begin{prop}\label{sec5:prop_bracket}
Let $(M,\omega,E,\nabla^E)$ be an $E$-valued $1$-plectic manifold. 
If an $E$-valued homotopy momentum section $\mu\in \Omega_A^1(M,E)$ is compatible with $A$, then it preserves the brackets up to the negative sign. 
Namely, it holds that 
\begin{equation}\label{sec5:eqninprop_bracket preserving}
\mu^{[\alpha,\beta]} = -\{\mu^\alpha,\,\mu^\beta\}
\end{equation}
for any section $\alpha,\beta$ of $A$. 
\end{prop}

\begin{rmk}
Consider the case where $E=\underline{\mathbb{R}}_M$ with the trivial connection $\mathrm{d}$ and $A$ is the action algebroid $A=\mathfrak{g}\ltimes M$. 
By \eqref{sec5:eqninprop_bracket preserving}, we have 
\[
 \alpha_M\mu^{\beta} = \mu^{[\alpha,\beta]}
\]
for any constant section $\alpha,\beta$ of $A$, where $\alpha_M=\rho(\alpha)$ is the infinitesimal generator of $\alpha$. 
From this, \eqref{sec5:eqninprop_bracket preserving} can be thought of as the operator like the coadjoint representation $\mathrm{ad}^*_{\alpha}\mu$. 
\end{rmk}
\medskip 

Denote by $\varGamma_{\mu}(E)$ the set of all sections of the form $\mu^\alpha\in \varGamma(E)$, i.e., 
\[
\varGamma_{\mu}(E) = \bigl\{\,\mu^\alpha\,|\,\alpha\in \varGamma(A)\,\bigr\}\subset \varGamma(E). 
\]
$\varGamma_{\mu}(E)$ is a subspace of ${\rm pHam}^0(M,E,\nabla^E)$.  From Proposition \ref{sec5:prop_bracket}, we find that $\omega$ defines a skew-symmetric bracket 
\begin{equation}\label{sec5:eqn_bracket2}
\{\cdot,\,\cdot\} : \varGamma_{\mu}(E) \times \varGamma_{\mu}(E) \longrightarrow \varGamma_{\mu}(E),\quad (\mu^\alpha,\,\mu^\beta)\longmapsto \omega\bigl(\rho(\alpha),\,\rho(\beta)\bigr) 
\end{equation} 
on $\varGamma_{\mu}(E)$. 

\begin{prop}\label{sec5:prop_bracket2}
$\varGamma_{\mu}(E)$ is a Lie algebra with respect to \eqref{sec5:eqn_bracket2}. 
\end{prop}
\begin{proof}
It suffices to show that the skew-symmetric bracket \eqref{sec5:eqn_bracket2} satisfies the Jacobi identity:
\[
\bigl\{\{\mu^\alpha,\,\mu^\beta\},\,\mu^\gamma\bigr\} + \bigl\{\{\mu^\beta,\,\mu^\gamma\},\,\mu^\alpha\bigr\} + \bigl\{\{\mu^\gamma,\,\mu^\alpha\},\,\mu^\beta\bigr\} = 0. 
\]
By Proposition \ref{sec5:prop_bracket}, it follows that
\begin{align*}
\omega\bigl(\rho([\alpha,\beta]),\,\rho(\gamma)\bigr)=\{\mu^{[\alpha,\beta]},\,\mu^{\gamma}\} = -\bigl\{\{\mu^\alpha,\,\mu^\beta\},\,\mu^\gamma\bigr\}, 
\end{align*}
and moreover, 
\begin{align*}
\nabla^E_{\rho(\gamma)}\{\mu^\alpha,\,\mu^\beta\} &=
-\nabla^E_{\rho(\gamma)}\mu^{[\alpha,\beta]}=-\bigl(\nabla^E\mu^{[\alpha,\beta]}\bigr)\bigl(\rho(\gamma)\bigr) = -\bigl(\imath_{\rho([\alpha,\beta])}\omega\bigr)(\rho(\gamma))\\
&= \bigl\{\{\mu^\alpha,\,\mu^\beta\},\,\mu^\gamma\bigr\}. 
\end{align*}
Accordingly, ${\rm d}_E^\nabla\omega$ is calculated as 
\begin{align*}
{\rm d}_E^\nabla\omega\bigl(\rho(\alpha),\rho(\beta),\rho(\gamma)\bigr)
&= \nabla^E_{\rho(\alpha)}\bigl(\omega(\rho(\beta),\rho(\gamma))\bigr) - \nabla^E_{\rho(\beta)}\bigl(\omega(\rho(\alpha),\rho(\gamma))\bigr) 
   + \nabla^E_{\rho(\gamma)}\bigl(\omega(\rho(\alpha),\rho(\beta))\bigr)\\
&\hspace{3.5em} -\omega\bigl(\rho([\alpha,\beta]),\,\rho(\gamma)\bigr) - \omega\bigl(\rho([\beta,\gamma]),\,\rho(\alpha)\bigr) + \omega\bigl(\rho([\alpha,\gamma]),\,\rho(\beta)\bigr)\\
&= \nabla^E_{\rho(\alpha)}\{\mu^\beta,\,\mu^\gamma\} - \nabla^E_{\rho(\beta)}\{\mu^\alpha,\,\mu^\gamma\} + \nabla^E_{\rho(\gamma)}\{\mu^\alpha,\,\mu^\beta\}\\
&\hspace{3.5em} + \bigl\{\{\mu^\alpha,\,\mu^\beta\},\,\mu^\gamma\bigr\} + \bigl\{\{\mu^\beta,\,\mu^\gamma\},\,\mu^\alpha\bigr\} + \bigl\{\{\mu^\gamma,\,\mu^\alpha\},\,\mu^\beta\bigr\}\\
&= \bigl\{\{\mu^\beta,\,\mu^\gamma\},\,\mu^\alpha\bigr\} + \bigl\{\{\mu^\gamma,\,\mu^\alpha\},\,\mu^\beta\bigr\} + \bigl\{\{\mu^\alpha,\,\mu^\beta\},\,\mu^\gamma\bigr\}\\
&\hspace{3.5em} + \bigl\{\{\mu^\alpha,\,\mu^\beta\},\,\mu^\gamma\bigr\} + \bigl\{\{\mu^\beta,\,\mu^\gamma\},\,\mu^\alpha\bigr\} + \bigl\{\{\mu^\gamma,\,\mu^\alpha\},\,\mu^\beta\bigr\}\\
&= 2\bigl(\bigl\{\{\mu^\alpha,\,\mu^\beta\},\,\mu^\gamma\bigr\} + \bigl\{\{\mu^\beta,\,\mu^\gamma\},\,\mu^\alpha\bigr\} + \bigl\{\{\mu^\gamma,\,\mu^\alpha\},\,\mu^\beta\bigr\}\bigr). 
\end{align*}
Since $\omega$ is ${\rm d}_E^\nabla$-closed, the assertion is proved. 
\end{proof}

\section{Reduction}
In the section, we shall establish the reduction theorem for an $E$-valued $1$-plectic manifold with an $E$-valued homotopy momentum section in the case where $E$ is a trivial bundle. 
\medskip 

Let $V$ be a finite dimensional real vector space and $(M,\omega,\underline{V}_M,\mathbf{d})$ an $\underline{V}_M$-valued 1-plectic manifold, 
where $\underline{V}_M$ denotes the trivial bundle $M\times V$ over $M$, and $\omega$ is a $V$-valued non-degenerate $2$-form on $M$ 
satisfying $\mathbf{d}\omega = \boldsymbol{0}_V$. The symbol $\boldsymbol{0}_V$ denotes the zero element in $V$. 
Note that $\mathbf{d}$ is given by 
\[
\mathbf{d}s = \sum_{i}({\rm d}s_i)\,\boldsymbol{v}_i
\] 
for any section $s:M\to V$ of $\underline{V}_M$ that has the form $s=\sum_i s_i\,\boldsymbol{v}_i,\,s_i\in C^{\infty}(M)$ with regard to a basis $\{\boldsymbol{v}_i\}$ of $V$. 
In addition, let $A$ be a Lie algebroid over $M$ endowed with a vector bundle connection $\nabla^A$, and 
$\mu \in \Omega_A^1(M,E)= \varGamma(A^*\otimes V)$ a $V$-valued homotopy momentum section. 



\subsection{The pseudogroup from the characteristic distribution of the Lie algebroid}

We let $\hat{0}$ be the zero section of $A^*\otimes V\cong \mathrm{Hom}(A,V)$ and denote by $M_{\mu}$ its preimage $\mu^{-1}(\hat{0})$ of $\hat{0}$. 
Assume that $M_\mu$ is an embedded submanifold of $M$ which is transversal to the image of the anchor map $\rho$ of $A$, i.e., 
\begin{equation}\label{sec5:cond_transversality}
 T_zM_{\mu} + {\rm Im}\,\rho_z = T_zM 
\end{equation} 
at each $z\in M_\mu$. Then, 
\[
A_{\mu}:= \bigl\{\,a\in A\,\bigm|\,\rho(a)\in TM_{\mu}\,\bigr\}
\]
is a Lie algebroid over $M_{\mu}$ by restricting $\rho$ to $A_\mu$ (which we denote by the same letter) and by defining the Lie bracket on $\varGamma(A_\mu)$ as 
\[
[\alpha_1,\,\alpha_2]:=\bigl.[\tilde{\alpha}_1,\,\tilde{\alpha}_2]\bigr|_{M_{\mu}}
\]
for $\alpha_1,\alpha_2\in \varGamma(A_{\mu})$, where $\tilde{\alpha}_1,\,\tilde{\alpha}_2$ are sections of $A$ that are extensions of them. 
\medskip 

The anchor map $\rho$ of $A_\mu$ yields a family of vector fields $\mathcal{C}_{\mu}$ on $M_\mu$ by 
\[
\mathcal{C}_{\mu} := {\rm im}\,\rho|_{M_{\mu}} = \bigl\{\,\rho(\alpha)\,|\,\alpha\in \varGamma(A_\mu)\,\bigr\}, 
\]
which gives rise to a singular distribution $\mathcal{D}_{\mu}$, called the characteristic distribution of $A_{\mu}$ 
\[
\mathcal{D}_{\mu} : M_{\mu} \ni z\longmapsto \mathcal{D}_{\mu}(z) := {\rm span}\bigl\{\,\rho(\alpha)_z \,|\,\alpha\in \varGamma(A_{\mu})\,\bigr\}\subset T_zM_{\mu}. 
\]
Denote by $\mathcal{E}_{\mu}$ the set of the local flows of a vector fields in $\mathcal{C}_{\mu}$, i.e., 
\[
\mathcal{E}_{\mu} = \{\,F_t^X\,|\, \text{the local flow of $X \in \mathcal{C}_{\mu}$} \,\}. 
\]
Then, we obtain the pseudogroup of transformations generated by it, 
\[
\mathcal{P}_{\mu} := \{{\rm id}_{M_{\mu}}\}\, \bigcup\, \Bigl\{\,F_{t_1}^{1}\circ\cdots\circ F_{t_k}^{k}\,\bigm|\,k\in\mathbb{N},\, 
 \text{$F_{t_j}^{j}\in \mathcal{E}_{\mu}$ or $(F_{t_j}^{j})^{-1}\in \mathcal{E}_{\mu}$} \,\Bigr\}.
\]
For each $z\in M_\mu$, we define the $\mathcal{P}_{\mu}$-orbit through $z$ to be the set 
\[
\mathcal{P}_{\mu}\cdot z := \bigl\{\,F_{\boldsymbol{t}}(z)\,|\,F_{\boldsymbol{t}}\in \mathcal{P}_{\mu},\, z\in {\rm Dom}\, F_{\boldsymbol{t}}\,\bigr\}, 
\]
where ${\rm Dom}\,F_{\boldsymbol{t}}$ denote the domain of the local diffeomorphism $F_{\boldsymbol{t}}$ in $\mathcal{P}_{\mu}$. 
Remark that $F_{\boldsymbol{t}}$ is expressed as $F_{\boldsymbol{t}}=F^{1}_{t_1}\circ\cdots\circ F^{k}_{t_k}$ for some $k$-tupple $\boldsymbol{t} = (t_1,\,\cdots,\,t_k)\in \mathbb{R}^k$. 
Two points $z$ and $z'$ in $M_\mu$ are said to be $\mathcal{P}_{\mu}$-equivalent if there exists an element $F_{\boldsymbol{t}}\in \mathcal{P}_{\mu}$ such that 
$z'=F_{\boldsymbol{t}}(z)$. The relation $\sim_{\mu}$ being $\mathcal{P}_{\mu}$-equivalent is an equivalence relation. 
The $\mathcal{P}_{\mu}$-orbit through $z$ coincides with the equivalence class $[z]_{\mu}$ of $z$. 
The orbit space, denoted by $\mathcal{M}_{\mu}:=M_\mu/\mathcal{P}_{\mu}$, is a topological space by the quotient topology. 
We use the notation $\pi_{\mu}:M_\mu\to \mathcal{M}_{\mu}$ for the canonical projection. 
\smallskip

The next result is due to \cite{Sorb73,Sacc74}. 

\begin{thm}[Stefan\cite{Sacc74} and Sussmann\cite{Sorb73}]\label{sec5;thm_stefan and sussmann}
Let $\mathscr{D}$ be a smooth singular distribution generated by a family $\mathscr{C}$ of smooth vector field. 
$\mathscr{P}_{\mathscr{C}}$ denotes the pseudogroup of transformations generated by the local flows of elements in $\mathscr{C}$. 
Then, the following conditions are equivalent{\rm :}
\begin{enumerate}[\quad \rm (1)]
 \item $\mathscr{D}$ is integrable and its maximal integral manifolds are the $\mathscr{P}_{\mathscr{C}}$-orbits. 
 \item At each point $x$, $\mathscr{D}(x)$ is the tangent space to the leaf, containing $x$, of the singular foliation associated to $\mathscr{D}$. 
 \item $\mathscr{D}$ is invariant with respect to $\mathscr{C}${\rm :} if $X\in \mathscr{C}$, then $(F_t^X)_x\mathscr{D}(x) = \mathscr{D}\bigl(F_t^X(x)\bigr)$ whenever the local flow 
  $F_t^X$ of $X$ is defined. 
\end{enumerate}
\end{thm}

The characteristic distribution of any Lie algebroid is integrable \cite{DZpoi05}. 
By Theorem \ref{sec5;thm_stefan and sussmann}, we have the following statement. 

\begin{cor}\label{sec5:cor_stefan and sussmann}
For each $F_{\boldsymbol{t}}\in \mathcal{P}_{\mu}$ and for each $z\in {\rm Dom}\,F_{\boldsymbol{t}}$, it holds that 
\[
({\rm d}F_{\boldsymbol{t}})_z(\mathcal{D}_{\mu}(z)) = \mathcal{D}_{\mu}\bigl(F_{\boldsymbol{t}}(z)\bigr). 
\]
Moreover, the tangent space $T_z(\mathcal{P}_{\mu}\cdot z)$ to the $\mathcal{P}_{\mu}$-orbit $\mathcal{P}_{\mu}\cdot z$ is written in the form
\[
 {\rm span}\Bigl\{\,({\rm d}F_{\boldsymbol{t}})_y\bigl(\rho_y(\alpha(y))\bigr)\, 
 \bigm|\,F_{\boldsymbol{t}}\in \mathcal{P}_{\mu},\,\alpha\in \varGamma(A_{\mu}),\,F_{\boldsymbol{t}}(y)= z\,\Bigr\}. 
\]
\end{cor}
For further discussion on the pseudogroup of transformations generated by local vector fields, we refer to Chapter 3 in \cite{ORmom05} for instance. 

\subsection{Vector bundle over the orbit space $\mathcal{M}_{\mu}$ with connection}

From the fact that the characteristic distribution $\mathcal{D}_{\mu}$ of $A_{\mu}$ is integrable, 
it follows that the orbit space $\mathcal{M}_{\mu}$ by $\mathcal{P}_{\mu}$ is a topological space. 
Assume that $\mathcal{M}_{\mu}$ has a smooth structure such that $\pi_{\mu}$ is smooth submersion, and 
consider the trivial bundle $\underline{V}_{\mathcal{M}_{\mu}}:= \mathcal{M}_{\mu}\times V$ over $\mathcal{M}_{\mu}$. 
Any section of $\underline{V}_{\mathcal{M}_{\mu}}$ is regarded as a $V$-valued function on $\mathcal{M}_{\mu}$. 

\begin{dfn}\label{sec5:dfn_the reduced function}
We say that a $V$-valued function $s\in C^{\infty}(M_{\mu}, V)$ is $\mathcal{P}_{\mu}$-invariant if it holds that $(s\circ F_{\boldsymbol{t}})(z) = s(z)$ 
for all $F_{\boldsymbol{t}}\in \mathcal{P}_{\mu}$ and $z\in \mathrm{Dom}\,F_{\boldsymbol{t}}$. 
We denote the set of $\mathcal{P}_{\mu}$-invariant $V$-valued function by $C^{\infty}(M_{\mu},V)^{\mathcal{P}_{\mu}}$, i.e., 
\[
 C^{\infty}(M_{\mu},V)^{\mathcal{P}_{\mu}} 
:= \bigl\{\,s\in C^{\infty}(M_{\mu}, V)\,|\, \forall F_{\boldsymbol{t}}\in \mathcal{P}_{\mu}~\text{and}~\forall z\in \mathrm{Dom}\,F_{\boldsymbol{t}}\, ;\, 
(s\circ F_{\boldsymbol{t}})(z) = s(z)\, \bigr\}. 
\]
Moreover, we define a function $\bar{s}\in C^{\infty}(\mathcal{M}_{\mu}, V)$ by 
\[
\bar{s}([z]_{\mu}) := s(z)~ ;\quad z\in M_{\mu}
\]
for $s\in C^{\infty}(M_{\mu},V)^{\mathcal{P}_{\mu}}$, and call it the reduced $V$-valued function. 
We denote by $C_{\rho}(\mathcal{M}_{\mu},V)$ the set of the reduced $V$-valued functions on $\mathcal{M}_{\mu}$. 
\end{dfn}

Similarly to the case of $V$-valued functions, we can define the notions of the $\mathcal{P}_{\mu}$-invariant function and the reduced $\mathbb{R}$-valued function by 
replacing $V$ with $\mathbb{R}$. We denote by $C^{\infty}(M_{\mu})^{\mathcal{P}_{\mu}}$ and $C_{\rho}(\mathcal{M}_{\mu})$, 
the set of $\mathcal{P}_{\mu}$-invariant functions on $M_{\mu}$ and that of the reduced functions on $\mathcal{M}_{\mu}$, respectively. 
Note that $s\in C^{\infty}(M_{\mu},V)^{\mathcal{P}_{\mu}}$ can be expressed in terms of $\mathcal{P}_{\mu}$-invariant function as 
\[
s=\sum_i s_i\,\boldsymbol{v}_i~ ;\quad s_i\in C^{\infty}(M_{\mu})^{\mathcal{P}_{\mu}}. 
\]

Let $\mathfrak{X}(\mathcal{M}_{\mu})$ denote the space of smooth vector fields on $\mathcal{M}_{\mu}$. 
Define a map $\widehat{\nabla}:\mathfrak{X}(\mathcal{M}_{\mu})\times C_{\rho}(\mathcal{M}_{\mu},V)\to C_{\rho}(\mathcal{M}_{\mu},V)$ by 
\begin{equation}\label{sec5:connection on reduced space}
(\widehat{\nabla}_{\bar{X}}\bar{s})\bigl([z]_{\mu}\bigr) := X_zs = (\mathbf{d}s)_z(X_z) \in V ~;\quad z\in M_{\mu}, 
\end{equation}
where $\bar{X}\in \mathfrak{X}(\mathcal{M}_{\mu}),\, X_z\in T_zM_{\mu}$ satisfying $\bar{X}_{[z]_{\mu}}=({\rm d}\pi_{\mu})_z(X_z)$, 
and where $\bar{s}$ is the reduced $V$-valued function for $s\in  C^{\infty}(M_{\mu},V)^{\mathcal{P}_{\mu}}$. 
The map \eqref{sec5:connection on reduced space} is shown to be a vector bundle connection on the trivial bundle $\underline{V}_{\mathcal{M}_{\mu}}$ in the following lemma. 

\begin{lem}\label{sec5:lem_connection on the reduced space}
The map $\widehat{\nabla}$ is well-defined, and satisfies that 
\begin{enumerate}[\quad \rm (1)]
 \item If $\bar{f}\in C_{\rho}(\mathcal{M}_{\mu})$, then 
  $(\widehat{\nabla}_{\bar{f}\bar{X}}\bar{s})\bigl([z]_{\mu}\bigr) = \bar{f}(\widehat{\nabla}_{\bar{X}}\bar{s})\bigl([z]_{\mu}\bigr)$. 
 \item If $\bar{f}\in C_{\rho}(\mathcal{M}_{\mu})$, then 
  $\bigl(\widehat{\nabla}_{\bar{X}}(\bar{f}\,\bar{s})\bigr)\bigl([z]_{\mu}\bigr) = (X_zf)\,s + \bar{f}(\widehat{\nabla}_{\bar{X}}\bar{s})\bigl([z]_{\mu}\bigr)$. 
\end{enumerate}
\end{lem}

\begin{proof}
Let $\bar{X}$ be any vector field on $\mathcal{M}_{\mu}$ and $z,\,z'$ any point in the same $\mathcal{P}_{\mu}$-orbit in $M_{\mu}$. 
Then, there exists an element $F_{\boldsymbol{t}}\in \mathcal{P}_{\mu}$ such that $z'=F_{\boldsymbol{t}}(z)$, and $({\rm d}\pi_{\mu})_z(X_z)=({\rm d}\pi_{\mu})_{z'}(X_{z'})$ 
for some $X_z\in T_zM_{\mu}$ and $X_{z'}\in T_{z'}M_{\mu}$. 
Since the projection $\pi_{\rho}$ is invariant under $F_{\boldsymbol{t}}$, we have 
\[
X_{z'} = ({\rm d}F_{\boldsymbol{t}})_z(X_z) + \sum_{i}h_i\bigl(F_{\boldsymbol{t}}(z)\bigr)\,\rho(\alpha_i)_{F_{\boldsymbol{t}}(z)},
\]
where $h_i\in C^\infty(M_{\mu})$ and $\alpha_i\in \varGamma(A_{\mu})$. If $s$ is a $\mathcal{P}_{\mu}$-invariant function, then, 
\[
\rho(\alpha_i)_zs = \rho(\alpha)_z(\bar{s}\circ\pi_{\mu})=(\mathbf{d}\bar{s})_{[z]_{\mu}}\left(({\rm d}\pi_{\mu})_z\bigl(\rho(\alpha_i)_z\bigr)\right) 
= \boldsymbol{0}_V 
\]
because every point on the integral curve of $\rho(\alpha_i)$ is projected to the same point of $\mathcal{M}_{\mu}$ by $\pi_{\mu}$. 
Therefore, 
\[
X_{z'}s = (\mathbf{d}s)_{F_{\boldsymbol{t}}(z)}\bigl(({\rm d}F_{\boldsymbol{t}})_z(X_z)\bigr) = (\mathbf{d}s)_z(X_z) = X_zs. 
\]
This shows that \eqref{sec5:connection on reduced space} is well-defined. The conditions (1) and (2) follows from the fact that both $f$ and $s$ are of $\mathcal{P}_{\mu}$-invariance. 
\end{proof}

The space of differential $k$-forms on $\mathcal{M}_{\mu}$, which we denote by $\Omega^k(\mathcal{M}_{\mu})$, becomes a $C_{\rho}(\mathcal{M}_{\mu})$-modules. 
The covariant exterior derivative $\mathbf{d}^{\widehat{\nabla}}$ of $\widehat{\nabla}$ is given by
\[
\mathbf{d}^{\widehat{\nabla}}=\mathbf{d}=\mathrm{d}\otimes 1:\Omega^k(\mathcal{M}_{\mu})\otimes V\to \Omega^{k+1}(\mathcal{M}_{\mu})\otimes V 
\] 
for $k=0,1,2,\cdots$, and satisfies that $\mathbf{d}^{\widehat{\nabla}}\circ\mathbf{d}^{\widehat{\nabla}}=0$. 


\subsection{Reduction with $V$-valued homotopy momentum section}

The non-degeneracy of the $V$-valued $1$-plectic form $\omega$ says that each of the induced linear maps 
\[
 \omega^{\flat}_x:T_xM\longrightarrow T_x^*M\otimes V,\quad \boldsymbol{u}\longmapsto \imath_{\boldsymbol{u}}\omega_x\, ;\quad x\in M. 
\] 
is injective. In other words, a vector $\boldsymbol{u}\in T_xM$ that satisfies $\omega_x(\boldsymbol{u},\,\boldsymbol{v})=\boldsymbol{0}_V$ for all $\boldsymbol{v}\in T_xM$ 
must be the zero vector. Extending this to general linear subspaces in $T_xM$ leads us to the notion of the $\omega$-orthogonality like the symplectic one in symplectic geometry. 

Let $W$ be a linear subspace of $T_xM$ at $x\in M$. We define the $\omega$-orthogonal subspace of $W$ as a linear subspace in $T_xM$ 
\[
W^\omega := \bigl\{\,\boldsymbol{u}\in T_xM\,|\,\forall \boldsymbol{v}\in W ;\, \omega_x(\boldsymbol{u},\boldsymbol{v})=\boldsymbol{0}_V\,\bigr\}. 
\]
In the case of $V=\mathbb{R}$, the $\omega$-orthogonal subspace of $W$ is nothing but the symplectic orthogonal subspace. In contrast with the symplectic case, it fails to hold 
that $W={\left(W^{\omega}\right)}^{\omega}$ in general. One can check easily the following proposition in the same manner as the case of $V=\mathbb{R}$. 

\begin{prop}
Let $x$ be a point in a $V$-valued $1$-plectic manifold $M$, and $W$ a linear subspace in $T_xM$. 
If $\dim V=1$, then it holds that $\dim T_xM=\dim W + \dim W^{\omega}$ and $W={\left(W^{\omega}\right)}^{\omega}$. 
\end{prop}

\begin{lem}\label{sec5:lem_subspaces}
Let $\mu\in \varGamma(A^*\otimes V)$ be a $V$-valued homotopy momentum section. 
The relationship $T_zM_{\mu}\subset \ker(\nabla^{A^*\otimes V}\mu)_z$ holds at each $z\in M_{\mu}$. 
Furthermore, $T_zM_{\mu}\subset \mathcal{D}_{\mu}(z)^\omega$ also holds. 
\end{lem}
\begin{proof}
Let $\boldsymbol{u}\in T_zM_{\mu}$. Then, there exists a smooth curve $\gamma$ in $M_{\mu}$ such that $\gamma(0)=z$ and $\frac{{\rm d}}{{\rm d}t}\gamma(0)=\boldsymbol{u}$. 
Then, the covariant derivative $\nabla^{A^*\otimes V}_{\dot{\gamma}(t)}\mu$ is given in terms of the covariant differentiation $\frac{\rm D}{\mathrm{d}t}$ along $\gamma$ by 
\[
 \nabla^{A^*\otimes V}_{\dot{\gamma}(t)}\mu =\frac{\rm D}{{\rm d}t}(\mu\circ\gamma)(t) ~;~ t\in \mathbb{R}. 
\]
If $\mu$ has the local form $\mu(t):=\mu(\gamma(t))=\sum_i\mu_i(t)\sigma_i(t)$ with regard to a smooth local frame $\{\sigma_i\}_i$ of $A^*\otimes V$, 
then $\frac{\rm D}{\mathrm{d}t}$ is expressed locally as 
\[
\frac{\rm D}{{\rm d}t}(\mu\circ\gamma)(t) = \sum_{i}\left(\frac{{\rm d}\mu^i}{{\rm d}t}(t)\,\sigma_i(t) 
+ \mu^i(t)\,\nabla^{A^*\otimes V}_{\dot{\gamma}(t)}\sigma_i(t)\right). 
\]
Since $(\mu\circ\gamma)(t)=\hat{0}$ for all $t$, it immediately follows that 
\[
\nabla^{A^*\otimes V}_{\boldsymbol{u}}\mu=\nabla^{A^*\otimes V}_{\dot{\gamma}(0)}\mu = \frac{\rm D}{{\rm d}t}(\mu\circ\gamma)(0) = \hat{0},
\]
which shows that $T_zM_{\mu}\subset \ker(\nabla^{A^*\otimes V}\mu)_z$. 
\smallskip 

Next, 
$\mu$ being a $V$-valued momentum section, we have 
\[
 (\nabla^{A^*\otimes V}_{\boldsymbol{v}}\mu)^\alpha = (\imath_\rho^1\omega)^\alpha(\boldsymbol{v}) = -\omega_z(\boldsymbol{v},\rho(\alpha)_z), 
\]
where $\boldsymbol{v}\in T_zM$, $\alpha\in\varGamma(A)$ and $z\in M_{\mu}$. 
From the equation, we immediately find that $\mathcal{D}_{\mu}(z)^\omega = \ker(\nabla^{A^*\otimes V}\mu)_z$. 
Combining it with the relation which we proved previously in the lemma, we see that $T_zM_{\mu}\subset \mathcal{D}_{\mu}(z)^\omega$. 
\end{proof}

\begin{dfn}
We say that the $V$-valued $1$-plectic form $\omega$ is $\mathcal{P}_{\mu}$-invariant if 
\begin{equation*}
 \forall F_{\boldsymbol{t}}\in \mathcal{P}_{\mu}~ \text{\rm and}~ \forall z\in {\rm Dom}\, F_{\boldsymbol{t}} ~;~ \left(F_{\boldsymbol{t}}^*\omega\right)_z = \omega_z. 
\end{equation*}
\end{dfn}
\medskip 

The reduction theorem for a vector-valued $1$-plectic manifold by a homotopy momentum section is stated as follows:

\begin{thm}\label{sec6:thm_main1}
Let $(M,\omega,\underline{V}_M,\mathbf{d})$ be a $V$-valued $1$-plectic manifold and $A$ a Lie algebroid over $M$ with a connection $\nabla^A$. 
Let $\mu\in \varGamma(A^*\otimes V)$ be a $V$-valued homotopy momentum section transversal to the characteristic distribution of $A$~{\rm (}see \eqref{sec5:cond_transversality}{\rm )}. 
Suppose that 
the preimage $M_{\mu}=\mu^{-1}(\hat{0})$ of the zero section $\hat{0}$ of $A^*\otimes V$ admits a smooth structure, and that
$\omega$ is $\mathcal{P}_{\mu}$-invariant. 
Additionally, assume that the orbit space $\mathcal{M}_{\mu}:=M_{\mu}/\mathcal{P}_{\mu}$ is a smooth manifold such that the canonical projection 
$\pi_{\mu}:M_\mu\to \mathcal{M}_{\mu}$ is smooth submersion. 

Then, $\mathcal{M}_{\mu}$ together with the trivial bundle $\underline{V}_{\mathcal{M}_{\mu}}$ with the connection $\widehat{\nabla}$ 
is a $V$-valued pre-$1$-plectic manifold whose $V$-valued {\em pre}-$1$-plectic form $\omega_{\mu}$ is uniquely characterized by 

\begin{equation}\label{sec5:eqn_relation_main thm}
\pi_{\mu}^*\omega_{\mu} = \imath_{\mu}^*\omega, 
\end{equation}
where $\imath_{\mu}:M_{\mu}\hookrightarrow M$ is the inclusion. 
\end{thm}

\begin{proof}

Define a $V$-valued 2-form $\omega_{\rm red}$ on $\mathcal{M}_{\mu}$ by 
\begin{equation}\label{sec5:eqn in thm_2-form}
(\omega_{\mu})_{[z]_{\mu}}(\bar{\boldsymbol{u}},\, \bar{\boldsymbol{v}}) := \omega_z(\boldsymbol{u},\,\boldsymbol{v})\in V, 
\end{equation} 
where $\boldsymbol{u},\, \boldsymbol{v}$ is any tangent vector of $M_{\mu}$ to $z$ and 
where $\bar{\boldsymbol{u}},\, \bar{\boldsymbol{v}}$ denote the equivalent classes of each of them in 
$T_{[z]_{\mu}}\mathcal{M}_{\mu}\cong T_zM_{\mu}/\mathcal{D}_{\mu}(z)$: $\bar{\boldsymbol{u}}=({\rm d}\pi_{\mu})_z(\boldsymbol{u}),\, \bar{\boldsymbol{v}}=({\rm d}\pi_{\mu})_z(\boldsymbol{v})$. 
To check that the expression \eqref{sec5:eqn in thm_2-form} is well-defined, we let $z'$ be a point in the $\mathcal{P}_{\mu}$-orbit through $z$, 
and let $\boldsymbol{u}',\,\boldsymbol{v}'\in T_{z'}M_{\mu}$ such that $({\rm d}\pi_{\mu})_{z'}(\boldsymbol{u}') = ({\rm d}\pi_{\mu})_{z}(\boldsymbol{u})$ and 
$({\rm d}\pi_{\mu})_{z'}(\boldsymbol{v}') = ({\rm d}\pi_{\mu})_{z}(\boldsymbol{v})$. 
By the same reasoning as in the proof of Lemma \ref{sec5:lem_connection on the reduced space}, it follows that 
\begin{align*}
\boldsymbol{u}' = ({\rm d}F_{\boldsymbol{t}})_z(\boldsymbol{u}) + \sum_{i}h_i\bigl(F_{\boldsymbol{t}}(z)\bigr)\,\rho(\alpha_i)_{F_{\boldsymbol{t}}(z)}, \quad
\boldsymbol{v}' = ({\rm d}F_{\boldsymbol{t}})_z(\boldsymbol{v}) + \sum_{j}k_j\bigl(F_{\boldsymbol{t}}(z)\bigr)\,\rho(\beta_j)_{F_{\boldsymbol{t}}(z)}, 
\end{align*}
where $F_{\boldsymbol{t}}\in \mathcal{P}_{\mu}$ such that $z'=F_{\boldsymbol{t}}(z)$, $h_i,\,k_j\in C^\infty(M_{\mu})$ and $\alpha_i,\,\beta_j\in \varGamma(A_{\mu})$. 
Using Lemma \ref{sec5:lem_subspaces}, we have 
\[
\omega_{z'}\bigl(({\rm d}F_{\boldsymbol{t}})_z(\boldsymbol{u}),\, \rho(\beta_j)_{z'}\bigr) = 
\omega_{z'}\bigl(({\rm d}F_{\boldsymbol{t}})_z(\boldsymbol{v}),\, \rho(\alpha_i)_{z'}\bigr) = \boldsymbol{0}_V 
\]
for each $i,j$. Thus, from the assumption that $\omega$ is $\mathcal{P}_{\mu}$-invariant, 
\begin{align}
\omega_{z'}(\boldsymbol{u}',\,\boldsymbol{v}') 
&= \omega_{z'}\bigl(({\rm d}F_{\boldsymbol{t}})_z(\boldsymbol{u}),\, ({\rm d}F_{\boldsymbol{t}})_z(\boldsymbol{v})\bigr) 
 + \sum_{i,j}h_i(z')k_j(z')\, \omega_{z'}\bigl(\rho(\alpha)_{z'},\,\rho(\beta)_{z'}\bigr) \notag \\
&= \omega_z(\boldsymbol{u},\,\boldsymbol{v}) 
 + \sum_{i,j}h_i\bigl(F_{\boldsymbol{t}}(z)\bigr)k_j\bigl(F_{\boldsymbol{t}}(z)\bigr)\, \{\mu^{\alpha},\mu^{\beta}\}\bigl(F_{\boldsymbol{t}}(z)\bigr). \label{sec5:eqn_proof of main thm}
\end{align}
Since $\mu^{\alpha}\equiv\boldsymbol{0}_V$ on $M_{\mu}$ for every $\alpha\in\varGamma(A_{\mu})$, we see that 
$\{\mu^{\alpha},\mu^{\beta}\}\bigl(F_{\boldsymbol{t}}(z)\bigr)=\boldsymbol{0}_V$. 
Substituting this to \eqref{sec5:eqn_proof of main thm}, we have $\omega_{z'}(\boldsymbol{u}',\,\boldsymbol{v}')=\omega_{z}(\boldsymbol{u},\,\boldsymbol{v})$, 
which shows that the 2-form $\omega_{\mu}$ is well-defined. 
In addition, the condition \eqref{sec5:eqn_relation_main thm} follows directly from \eqref{sec5:eqn in thm_2-form}. 

Lastly, we verify that the $V$-valued 2-form $\omega_{\mu}$ is closed with respect to $\mathbf{d}^{\widehat{\nabla}}$. 
Note that $\mathbf{d}^{\widehat{\nabla}}$ commutes with the pullback map $\pi_{\mu}^*$ because $\mathbf{d}^{\widehat{\nabla}}=\mathrm{d}\otimes 1$. 
From \eqref{sec5:eqn_relation_main thm} and $\mathbf{d}\omega=\boldsymbol{0}_V$, it follows that 
\[
\pi_{\mu}^*\mathbf{d}^{\widehat{\nabla}}\omega_{\mu} = \mathbf{d}\imath_{\mu}^*\omega = \imath_{\mu}\mathbf{d}\omega = \boldsymbol{0}_V. 
\]
This implies that the $V$-valued 2-form $\omega_{\mu}$ is closed with respect to $\mathbf{d}^{\widehat{\nabla}}$. This completes the proof. 
\end{proof}
\subsection{Reduction with the $V$-valued homotopy momentum section compatible with Lie algebroid}

In the subsection, we do {\em not} require that the $V$-valued homotopy momentum section $\mu$ is transversal to the anchor map $\rho$ of $A$~
(see \eqref{sec5:cond_transversality}). 
Instead of that, we assume that $\mu$ is compatible with $A$~(see Definition \ref{sec5:dfn_compatibility with A}). 
We let $M_{\mu}$ be the preimage of the zero section $\hat{0}$ of $A^*\otimes V$, and suppose that it is an embedded submanifold of $M$. 
\medskip 

$\mathcal{C}_{\rho}$ denotes the family of vector field on $M$ consisting of all the images of 
the anchor map $\rho$ of $A$, and moreover, $\mathcal{D}_{\rho}$ denotes the characteristic distribution on $M$ associated to $\mathcal{C}_{\rho}$. That is, 
\[
\mathcal{C}_{\rho} := {\rm im}\,\rho = \bigl\{\,\rho(\alpha)\,|\,\alpha\in \varGamma(A)\,\bigr\},
\]
and 
\[
\mathcal{D}_{\rho}(x) := {\rm span}\bigl\{\,\rho(\alpha)_x\,|\,\alpha\in \varGamma(A)\,\bigr\}\subset T_xM~;\quad x\in M. 
\]
Following the discussion in the subsection 6.1, we can obtain the pseudogroup of transformations generated by the local flows of vector fields in $\mathcal{C}_{\rho}$, 
\[
\mathcal{P}_{\rho} := \{{\rm id}_{M}\}\, \bigcup\, \Bigl\{\,G_{t_1}^{1}\circ\cdots\circ G_{t_k}^{k}\,\bigm|\,k\in\mathbb{N},\, 
 \text{$G_{t_j}^{j}\in \mathcal{E}_{\rho}$ or $(G_{t_j}^{j})^{-1}\in \mathcal{E}_{\rho}$} \,\Bigr\}, 
\]
where $\mathcal{E}_{\rho}$ is the set of the local flows of a vector fields in $\mathcal{C}_{\mu}$: 
$\mathcal{E}_{\rho} = \{\,G_t^X\,|\, \text{the local flow of $X \in \mathcal{C}_{\rho}$} \,\}$. 
The $\mathcal{P}_{\rho}$-orbits is also defined in the same way as the case of $\mathcal{P}_{\mu}$. 
We write $G_{\boldsymbol{t}}$ for the elements in $\mathcal{P}_{\rho}$ that has the form 
$G_{\boldsymbol{t}}=G^{1}_{t_1}\circ\cdots\circ G^{k}_{t_k}$, where $\boldsymbol{t} = (t_1,\,\cdots,\,t_k) \in \mathbb{R}^k$. 
We say that two points $x$ and $x$ in $M$ are $\mathcal{P}_{\rho}$-equivalent if they are in the same $\mathcal{P}_{\rho}$-orbit, and denote by $[x]_{\rho}$ the equivalence class 
including $x$. In other words, $x$ and $x'$ are $\mathcal{P}_{\rho}$-equivalent if and only if there exists an element $G_{\boldsymbol{t}}$ in $\mathcal{P}_{\rho}$ satisfying 
$x'=G_{\boldsymbol{t}}(x)$. 
\medskip 

Since the singular distribution $\mathcal{D}_{\rho}$ is integrable, $\mathcal{D}_{\rho}$ is invariant with respect to $\mathcal{P}_{\rho}$: 
\begin{equation}\label{sec6:eqn_distribution_invariant}
({\rm d}G_{\boldsymbol{t}})_x(\mathcal{D}_{\rho}(x)) = \mathcal{D}_{\rho}\bigl(G_{\boldsymbol{t}}(x)\bigr)~;\quad 
G_{\boldsymbol{t}}\in \mathcal{P}_{\rho},\quad x\in {\rm Dom}\,G_{\boldsymbol{t}}. 
\end{equation}
Furthermore, Theorem \ref{sec5;thm_stefan and sussmann} guarantees that the maximal integral manifold through $x\in M$ of $\mathcal{D}_{\rho}$ is the leaf $\mathcal{L}_{\rho}(x)$, 
containing $x$, of the singular foliation associated to $\mathcal{D}_{\rho}$, and coincides with the $\mathcal{P}_{\rho}$-orbit $\mathcal{P}_{\rho}\cdot x$. 
Furthermore, 
\[
 T_x(\mathcal{P}_{\rho}\cdot x) = 
 {\rm span}\Bigl\{\,({\rm d}G_{\boldsymbol{t}})_y\bigl(\rho(\alpha)_y\bigr)\, 
 \bigm|\,G_{\boldsymbol{t}}\in \mathcal{P}_{\rho},\, \alpha\in \varGamma(A),\,G_{\boldsymbol{t}}(y)= x\,\Bigr\}. 
\]
\smallskip 

Define the subsets $ \mathcal{P}^{0}_{\rho}\subset \mathcal{P}_{\rho}$ and $\mathcal{P}^{0}_{\rho}\cdot x\subset \mathcal{P}_{\rho}\cdot x$ as 
\[
 \mathcal{P}^{0}_{\rho} := \{\, G_{\boldsymbol{t}}\in \mathcal{P}_{\rho}\,|\, \mu\circ G_{\boldsymbol{t}} = \hat{0}|_{\mathrm{Dom}\, G_{\boldsymbol{t}}}\,\}
\]
and 
\[
 \mathcal{P}^{0}_{\rho}\cdot x := \{\, G_{\boldsymbol{t}}(x)\,|\, 
 \text{$G_{\boldsymbol{t}}\in \mathcal{P}^0_{\rho},\, \mathrm{Dom}\, G_{\boldsymbol{t}}\ni x$ and $\mu(G_{\boldsymbol{t}}(x)) = 0\in A_x^*\otimes V$}\, \}~;~ x\in M
\]
respectively. The following proposition is easily shown. 

\begin{prop}\label{sec6:prop_invariance}
The submanifold $M_{\mu}$ in $M$ is $\mathcal{P}^0_{\rho}$-invariant, i.e., $\mathcal{P}^0_{\rho}\cdot z\subset M_{\mu}$ if $z\in M_{\mu}$. 
Furthermore, $\mathcal{P}^0_{\rho}\cdot z = \mathcal{L}_{\rho}(z) \cap M_{\mu}$ holds for any $z\in M_{\mu}$. 
\end{prop}

The proposition says that each $\mathcal{P}^0_{\rho}$-orbit through a point in $M_{\mu}$ are entirely included in $M_{\mu}$. 
Note that the characteristic distribution $\mathcal{D}_{\rho}$ is integrable. For every point $x$ in $M_{\mu}$, there exists a unique connected integral manifold $\mathcal{L}_{\rho}(x)$, 
the leaf of $\mathcal{D}_{\rho}$. 
We denote by $\mathcal{M}_{\rho}:=M_{\mu}/\mathcal{P}^0_{\rho}$ the orbit space by $\mathcal{P}^0_{\rho}$. 
We suppose that $\mathcal{M}_{\rho}$ is a smooth manifold such that the canonical projection $\pi_{\rho}:M_{\mu}\to \mathcal{M}_{\rho}$ is a surjective submersion. 
\medskip 

In the same manner as the subsection 6.2, the notions of $\mathcal{P}^0_{\rho}$-invariant functions and the reduced functions are again introduced. 
Namely, a $\mathcal{P}^0_{\rho}$-invariant $V$-valued function is a $V$-valued function $s\in C^{\infty}(M_{\mu}, V)$ satisfying $(s\circ G_{\boldsymbol{t}})(z) = s(z)$ 
for all elements $G_{\boldsymbol{t}}$ in $\mathcal{P}^0_{\rho}$ and $z\in \mathrm{Dom}\,G_{\boldsymbol{t}}$. 
We denote the set of $\mathcal{P}^0_{\rho}$-invariant $V$-valued function by $C^{\infty}(M_{\mu},V)^{\mathcal{P}^0_{\rho}}$. 
The reduced $V$-valued function is defined to be a function $\bar{s}\in C^{\infty}(\mathcal{M}_{\rho}, V)$ by 
$\bar{s}([z]_{\rho}) := s(z),\, z\in M_{\mu}$ for $s\in C^{\infty}(M_{\mu},V)^{\mathcal{P}^0_{\rho}}$. 

Define a map $\widehat{\nabla}$ for $\bar{X}\in \mathfrak{X}(\mathcal{M}_{\rho})$ and 
the reduced $V$-valued function $\bar{s}$ for $s\in  C^{\infty}(M_{\mu},V)^{\mathcal{P}^0_{\rho}}$ in the same method as \eqref{sec5:connection on reduced space}:
\begin{equation*}
(\widehat{\nabla}_{\bar{X}}\bar{s})\bigl([z]_{\rho}\bigr) := X_zs = (\mathbf{d}s)_z(X_z) \in V ~;\quad z\in M_{\mu}, 
\end{equation*}
where $X_z$ is the tangent vector satisfying $\bar{X}_{[z]_{\rho}}=({\rm d}\pi_{\rho})_z(X_z)$. 
By the same proof of Lemma \ref{sec5:lem_connection on the reduced space}, we have the following proposition:

\begin{prop}
The map $\widehat{\nabla}$ defines a vector bundle connection on the trivial bundle $\underline{V}_{\mathcal{M}_{\rho}}$. 
\end{prop}
Similarly to Theorem \ref{sec6:thm_main1}, the orbit space $\mathcal{M}_{\rho}$ is also proven to be a $V$-valued pre-1-plectic manifold that is stated as follows:

\begin{thm}\label{sec6:thm_main2}
Let $(M,\omega,\underline{V}_M,\mathbf{d})$ be a $V$-valued $1$-plectic manifold and $A$ a Lie algebroid over $M$ with a connection $\nabla^A$. 
Let $\mu\in \varGamma(A^*\otimes V)$ be a $V$-valued homotopy momentum section compatible with $A$. 
Suppose that $M_{\mu}=\mu^{-1}(\hat{0})$ is a smooth manifold, and that $\omega$ is $\mathcal{P}^0_{\rho}$-invariant, i.e., 
\begin{equation*}
 \forall G_{\boldsymbol{t}}\in \mathcal{P}^0_{\rho}~ \text{\rm and}~ \forall z\in {\rm Dom}\, G_{\boldsymbol{t}} ~;~ \left(G_{\boldsymbol{t}}^*\omega\right)_z = \omega_z. 
\end{equation*}
Furthermore, assume that the orbit space $\mathcal{M}_{\rho}:=M_{\mu}/\mathcal{P}^0_{\rho}$ is a smooth manifold such that the canonical projection 
$\pi_{\rho}:M_\mu\to \mathcal{M}_{\rho}$ is smooth submersion. 

Then, $\mathcal{M}_{\rho}$ together with the trivial bundle $\underline{V}_{\mathcal{M}_{\rho}}$ with the connection $\widehat{\nabla}$ 
is a $V$-valued pre-$1$-plectic manifold whose $V$-valued {\em pre}-$1$-plectic form $\omega_{\rho}$ is uniquely characterized by 

\begin{equation}\label{sec6:eqninproof:thm_main2}
\pi_{\rho}^*\omega_{\rho} = \imath_{\mu}^*\omega. 
\end{equation}
\end{thm}

\begin{proof}
This is shown by the same manner as the proof of Theorem \ref{sec6:thm_main1}. 
\end{proof}

\begin{ex}\label{sec6:ex_symplectic reduction}
Every symplectic manifold $(M,\,\omega)$ is an $\underline{\mathbb{R}}_M$-valued $1$-plectic manifold with the trivial connection $\nabla^E=\mathrm{d}$ 
{\rm (}Example \ref{sec2:ex_symplectic}{\rm )}.  
Assume that $(M,\omega)$ admits a Hamiltonian $G$-action $\Phi:G\times M\to M$ by a compact Lie group $G$.  
The momentum map $J:M\to \mathfrak{g}^*$ is an $\mathfrak{g}$-valued homotopy momentum section with respect to the action algebroid $A=\mathfrak{g}\ltimes M$~
{\rm (}Example \ref{sec3:example_momentum maps for symplectic manifolds}{\rm )}. 
Note that $J$ is equivariant with respect to the $G$-action, that is, $J(\Phi_g(x))=\mathrm{Ad}^*_gJ(x)$ is satisfied for all $x\in M$. 
Additionally, remark that $J$ is compatible with $A$, confined to the constant sections. 

Suppose that the zero element $\boldsymbol{0}\in\mathfrak{g}$ is a regular value of $J$ and set $M_{J}=J^{-1}(\boldsymbol{0})$. 
Since $J$ is equivariant with respect to the $G$-action, one finds that each $G$-orbit $G\cdot z$ is included in $M_J$ if $z\in M_J$. 
Consider the family of vector fields on $M$ consisting of all the infinitesimal generators, 
$\mathcal{C}_{\rho} = \{\,\xi_M\,|\, \xi\in \mathfrak{g}\,\}$. 
Then, the characteristic distribution $\mathcal{D}_{\rho}$ associated to $\mathcal{C}_{\rho}$ is given by 
\[
\mathcal{D}_{\rho}(x)=\mathrm{span}\{(\xi_M)_x\,|\,\xi\in\mathfrak{g}\}=T_x(G\cdot x)~;~ x\in M. 
\]
$G$ acting canonically on $M$, the symplectic form $\omega$ is $\mathcal{P}_{\rho}$-invariant, where $\mathcal{P}_{\rho}$ is the pseudogroup generated by $\mathcal{C}_{\rho}$. 
By Theorem \ref{sec5;thm_stefan and sussmann}, the maximal integral manifolds of $\mathcal{D}_{\rho}$ are the $G$-orbits. Therefore, from Theorem \ref{sec6:thm_main2} 
we obtain a presymplectic manifold $M_{J}/\mathcal{P}_{\rho}=M_{J}/G$ by \eqref{sec6:eqninproof:thm_main2}. In this case, it is shown that $\omega_{\rho}$ is nondegenerate. 
$M_{J}/G$ is none other than the symplectic manifold from the Marsden-Weinstein-Meyer reduction \cite{Msym73, MWred74}. 
\end{ex}

\begin{ex}
Consider a hyper-K\"{a}hler manifold $M$ equipped with a hyper-hamiltonian $G$-action. Recall that $M$ is an $\underline{\mathbb{R}^3}_M$-valued $1$-plectic manifold 
by $\omega^{(3)}=\sum_i^3\omega_i\otimes \boldsymbol{e}_i$ {\rm (}Example \ref{sec2:ex_family of symplectic}{\rm )} 
and the momentum map for the action $\mu=\sum_i^3\mu_i\otimes \boldsymbol{e}_i$ is thought of as an $\underline{\mathbb{R}^3}_M$-valued 
homotopy momentum section {\rm (}see Proposition \ref{sec4:prop_hyperkahler momentum maps}{\rm )}. 

Put $M_{\mu}:= \mu_1^{-1}(\boldsymbol{0}) \cap \mu_2^{-1}(\boldsymbol{0})\cap \mu_2^{-1}(\boldsymbol{0})$ for $\boldsymbol{0}\in \mathfrak{g}$. 
From the condition that each $\mu_i$ is $G$-equivariant, it follows that $\mu_1(\Phi_h(z))=\mu_2(\Phi_h(z))=\mu_3(\Phi_h(z))=\boldsymbol{0}$ for any $h\in G$ and $z\in M_{\mu}$. 
Namely, each $G$-orbit through a point in $M_{\mu}$ is included in $M_{\mu}$. 
The characteristic distribution $\mathcal{D}_{\rho}$ from the action Lie algebroid $A=\mathfrak{g}\ltimes G$ is given by the same as Example \ref{sec6:ex_symplectic reduction}. 
Since each symplectic form $\omega_i$ is $G$-invariant, the $\underline{\mathbb{R}^3}_M$-valued $1$-plectic form $\omega^{(3)}$ is invariant under the pseudogroup $\mathcal{P}_\rho$ 
associated to $\mathcal{D}_{\rho}$. 
By Theorem \ref{sec6:thm_main2}, the quotient manifold $\mathcal{M}_{\rho}:=M_{\mu}/\mathcal{P}_{\rho}=M_{\mu}/G$ is a $\underline{\mathbb{R}^3}_M$-valued pre-$1$-plectic 
manifold by \eqref{sec6:eqninproof:thm_main2}. 
In fact, $\mathcal{M}_{\rho}$ is proven to be a hyper-K\"{a}hler manifold again \cite{HKLR87}. 
\end{ex}

\begin{ex}
Let $M=\mathbb{T}^4$ be the $4$-torus, and $V=\mathbb{R}^3$ with the trivial connection. We denote by $(\theta^0,\theta^1,\theta^2,\theta^3)$ angle coordinates on $\mathbb{T}^4$. 
A circle $S^1$ naturally acts on $\mathbb{T}^4$ by the translation on the first component{\rm :} 
\[
\Phi_{\theta}(\theta_0,\theta_1,\theta_2,\theta_3)=(\theta + \theta^0,\, \theta^1,\, \theta^2,\, \theta^3).
\] 
Define a $\mathbb{R}^3$-valued $1$-plectic form $\omega$ and an $\mathbb{R}^3$-valued function $\mu:\mathbb{T}^4\to \mathrm{Lie}(S^1)\otimes \mathbb{R}^3$ by
\[
\omega = (\mathrm{d}\theta^0\wedge \mathrm{d}\theta^1)\,\boldsymbol{e}_1 + (\mathrm{d}\theta^1\wedge \mathrm{d}\theta^2)\,\boldsymbol{e}_2 
 + (\mathrm{d}\theta^1\wedge \mathrm{d}\theta^3)\,\boldsymbol{e}_3, 
\]
and 
\[
\mu^{\alpha}(p) = \langle \mu(p),\,\alpha\rangle = \theta^1\,\boldsymbol{e}_1;\quad p=(\theta^0,\theta^1, \theta^2, \theta^3)\in \mathbb{T}^4,\, \alpha\in \mathrm{Lie}(S^1), 
\]
respectively. Here $\{\boldsymbol{e}_i\}_{i=1}^3$ denotes the standard basis of $\mathbb{R}^3$. $\omega$ is invariant under the $S^1$-action. 
When if we restrict our attention to the constant section of 
the action algebroid $A=\mathrm{Lie}(S^1)\ltimes \mathbb{T}^4$, $\mu$ is an $\mathbb{R}^3$-valued homotopy momentum section compatible with $A$. 
The pseudo-Hamiltonian vector field corresponding to $\mu^{\alpha}$ is given by $\rho(\alpha)=\partial/\partial \theta^0$. 

Then, the preimage $M_{\mu}$ is a manifold expressed locally in the form $\{(\theta^0,0,\theta^2,\theta^3)\}$. 
By Theorem \ref{sec6:thm_main2}, the reduced space $\mathcal{M}_{\rho}$ is a manifold locally diffeomorphic to a $2$-torus $\mathbb{T}^2$, and 
the reduced $2$-form $\omega_{\rho}$ is zero. Namely, $\mathcal{M}_{\rho}$ is a $\mathbb{R}^3$-valued pre $1$-plectic manifold. 

\end{ex}
\section{Conclusions}

In the paper, we have introduced geometric objects in order to understand Hamiltonian symmetries in a single framework --- 
a bundle-valued (pre-)$n$-plectic structures and a bundle-valued homotopy momentum section (BHMS, for short). 
In contrast to conventional Cartan calculus, 
the curvatures of the vector bundle connections influence substantially the Cartan formulas for bundle-valued $n$-plectic manifolds. 
Both hyper-K\"{a}hler manifold and quaternionic K\"{a}hler manifold are regarded as bundle-valued $1$-plectic manifolds. 
It turns out that a hyper-K\"{a}hler momentum map is naturally a BHMS. 
Subsequently, a quaternionic K\"{a}hler momentum map is a BHMS under the condition that it is a Lie anti-homomorphism. 
The authors expect a BHMS to be the candidate for the framework integrating various momentum map theories. 

Furthermore, the study also describes new generalization of the Marsden-Weinstein-Meyer reduction for vector-valued $1$-plectic manifolds with Lie algebroid symmetry. 
We have constructed two kinds of reductions for obtaining a vector-valued pre-$1$-plectic manifold in Theorem \ref{sec6:thm_main1} and \ref{sec6:thm_main2}. 
Those are applicable to the case only of vector-valued $1$-plectic manifolds. We need to address the reduction for the case of
a general bundle-valued $n$-plectic manifold with BHMS. 

In addition, the study needs further investigation: the first is to find the condition for the reduced 2-forms $\omega_{\mu},\,\omega_{\rho}$ on the reduced spaces 
$\mathcal{M}_{\mu},\,\mathcal{M}_{\rho}$ to be non-degenerate. 
It might be easy to show that those 2-forms are non-degenerate if the vector bundle $E$ has the rank equal to $1$. 
However, the same manner is not applicable for the case of the rank greater than $1$. 
It would be important to see what kind of condition is needed for non-degeneracy. 

The second is to establish the reduction theory involving Lie groupoid action. Assuming that the Lie algebroid $A$ is integrable and acts on a bundle-valued $n$-plectic manifold, 
we obtain a Lie groupoid $\mathcal{G}$ associated to $A$ and $\mathcal{G}$-action integrating the action of $A$. 
It is worthwhile mathematically to discuss the reduction by $\mathcal{G}$-action. 

The third is to define the notion of the Hamiltonian Lie algebroid \cite{BWham18} for a Lie algebroid with a BHMS, and to 
compare their differences between a Lie algebroid in the case of BHMS for $n=1$ and that in the case of momentum section. 
As mentioned in Example \ref{sec3:example_momentum sections}, every momentum section $\mu$ is BHMS even if the Lie algebroid $A$ is not necessarily presymplectically anchored, i.e., 
$R_{\nabla}^A\mu = 0$. 
It would be interesting to find the mathematical meaning of the condition that a Lie algebroid is presymplectically anchored in the case of BHMS. 

The authors hope to address some of those problems in the future, 
and expect this paper provides insights into the study of multisymplectic geometry and the Hamiltonian symmetries.


\subsection*{Acknowledgment}
This work was supported by the research promotion program for acquiring grants in-aid for JSPS KAKENHI Grant Number 22K03323. 

\end{document}